\documentclass[a4paper]{article}

\usepackage{amssymb}
\usepackage{latexsym}
\usepackage{amsfonts}
\usepackage{mathrsfs}
\usepackage{amscd}
\usepackage{amsmath}
\usepackage{amssymb}
\usepackage{wasysym}
\usepackage{bm}
\newtheorem{Th}{Theorem}
\newtheorem{Prop}{Proposition}[section]
\newtheorem{Co}{Corollary}
\newtheorem{Lm}{Lemma}[section]
\newtheorem{Lma}{Lemma}[section]
\newtheorem{Dfi}{Definition}
\newtheorem{Rm}{Remark}

\newcommand{\be}{\begin{equation}}
\newcommand{\ee}{\end{equation}}
\newcommand{\bes}{\begin{equation*}}
\newcommand{\ees}{\end{equation*}}
\newcommand{\R}{\mathbb{R}}
\newcommand{\N}{\mathbb{N}}
\newcommand{\C}{\mathbb{C}}

\newcommand\res{\mathop{\hbox{\vrule height 7pt width .5pt depth 0pt
\vrule height .5pt width 6pt depth 0pt}}\nolimits}

\def\theequation{\thesection.\arabic{equation}}
\def\theTh{\Roman{section}.\arabic{Th}}
\def\theProp{\Roman{section}.\arabic{Prop}}
\def\theCo{\Roman{section}.\arabic{Co}}
\def\theLm{\Roman{section}.\arabic{Lma}}

\def\theRm{\Roman{section}.\arabic{Rm}}
\newcommand{\reset}{\setcounter{equation}{0}\setcounter{Th}{0}\setcounter{Prop}{0}\setcounter{Co}{0}\setcounter{Lma}{0}\setcounter{Rm}{0}}

\def\lf{\left}
\def\rg{\right}

\def\al{a}
\def\la{\lambda}
\def\eps{\varepsilon}

\def\Om{\Omega}

\def\p{\partial}
\def\pro{\pi_{\vec{n}}}
\def\bn{\vec{n}}

\def\bh{\vec{h}}
\def\bex{\bAe_1}
\def\bey{\bAe_2}
\def\bez{\bAe_z}

\def\bei{\bAe_i}
\def\bej{\bAe_j}

\def\bezb{\bAe_{\bar{z}}}
\def\px{\partial_{x_1}}
\def\py{\partial_{x_2}}
\def\pj{\partial_{x_j}}

\def\pz{\partial_{z}}
\def\pzb{\partial_{\bar{z}}}

\def\bAe{\vec{e}}
\def\bH{\vec{H}}

\def\bC{\vec{C}}
\def\bA{\vec{A}}
\def\bAe{\vec{e}}

\def\bF{\vec{F}}
\def\bJ{\vec{J}}

\def\bL{\vec{L}}
\def\bG{\vec{G}}
\def\bR{\vec{R}}
\def\bE{\vec{E}}

\def\bX{\vec{X}}

\def\bX{\vec{X}}
\def\bPe{\vec{P}}

\def\bB{\vec{B}}

\def\bp{\vec{\Phi}}
\def\bP{\vec{\Phi}}

\def\bT{\vec{T}}

\def\bQ{\vec{Q}}

\def\bul{\bullet}
\def\di{{D}^2}

\def\res{\mathop{\hbox{\vrule height 7pt width .5pt 
depth 0pt\vrule height .5pt width 6pt depth 0pt}}\nolimits}

\begin{document}


\reset

\title{Analysis of Constrained Willmore Surfaces}
\author{Yann Bernard\footnote{Mathematisches Institut, Albert-Ludwigs-Universit\"at, 79104 Freiburg, Germany. {\it Supported by the DFG Collaborative Research Center SFB/Transregio 71 (Project B3)}.}}
\maketitle
\noindent
$\textbf{Abstract:}$ {\it This paper studies the regularity of constrained Willmore immersions into $\R^{m\ge3}$ locally around both ``regular" points and around branch points, where the immersive nature of the map degenerates. We develop  local asymptotic expansions for the immersion, its first, and its second derivatives, given in terms of residues which are computed as circulation integrals. We deduce explicit ``point removability" conditions ensuring that the immersion is smooth. Our results apply in particular to Willmore immersions and to parallel mean curvature immersions in any codimension.}

\section{Introduction}

\subsection{Preliminaries}

Let $\vec{\Phi}$ be an immersion from a closed abstract two-dimensional manifold $\Sigma$ into ${\R}^{m\ge3}$. We denote by $g:=\vec{\Phi}^\ast g_{{\R}^m}$ the pull back by $\vec{\Phi}$ of the flat  canonical metric $g_{{\R}^m}$ of ${\R}^m$, also called the {\it first fundamental form of $\vec{\Phi}$}, and we let $d\text{vol}_g$ be its associated {\it volume form}. The {\it Gauss map} of the immersion $\vec{\Phi}$ is the map taking values in the Grassmannian of oriented $(m-2)$-planes in ${\R}^m$ given by
\[
\vec{n}\;:=\;\star\,\frac{\p_{x_1}\vec{\Phi}\wedge\p_{x_2}\vec{\Phi}}{|\p_{x_1}\vec{\Phi}\wedge\p_{x_2}\vec{\Phi}|}\:,
\] 
where $\star$ is the usual Hodge star operator in the Euclidean metric, and $\{x_1,x_2\}$ are local coordinates on the surface $\Sigma$. \\
Denoting by $\pi_{\vec{n}}$ the orthonormal projection of vectors in ${\R}^m$ onto the $(m-2)$-plane given by $\vec{n}$, the {\it second fundamental form} may be expressed as\footnote{In order to define $d^2\vec{\Phi}(X,Y)$ one has to extend locally around $T_p\Sigma$ the vector-fields $X$ and $Y$. It is not difficult to
check that $\pi_{\vec{n}}\,d^2\vec{\Phi}(X,Y)$ is independent of this extension.}

\[
\vec{\mathbb I}_p(X,Y)\;:=\;\pi_{\vec{n}}\,d^2\vec{\Phi}(X,Y)\qquad\forall\:\: X,Y\in T_p\Sigma\:.
\]The {\it mean curvature vector} of the immersion at the point $p\in\Sigma$ is 
\[
\vec{H}\;:=\;\frac{1}{2}\,\text{Tr}_g(\vec{\mathbb I}_p)\;=\;\frac{1}{2}\,\lf[\vec{\mathbb I}_p(\bAe_1,\bAe_1)+\vec{\mathbb I}_p(\bAe_2,\bAe_2)\rg]\:,
\]
where $\{\bAe_1,\bAe_2\}$ is an orthonormal basis of $T_p\Sigma$ for the metric $g$. 

\medskip

In the present paper, we study the functional
\[
W(\vec{\Phi})\;:=\;\int_\Sigma|\vec{H}|^2\,d\text{vol}_{g}\:,
\]
called {\it Willmore energy}. It has been extensively studied in the literature, due to its relevance to various areas of science. We refer the reader to \cite{Ri3} and the references therein for more extensive information on the properties and applications of the Willmore energy. \\
The Gauss-Bonnet theorem and Gauss equation imply that 
\bes
W\big(\bP\big)\;=\;\dfrac{1}{4}\int_{\Sigma}\,\big|\vec{\mathbb{I}}\big|_g^2\,d\text{vol}_g+\,\pi\chi(\Sigma)\;=\;\dfrac{1}{4}\int_{\Sigma}\,\big|d\bn|_g^2\,d\text{vol}_g+\,\pi\chi(\Sigma)\:,
\ees
where $\chi(\Sigma)$ is the Euler characteristic of $\Sigma$, which is a topological invariant for a closed surface. From the variational point of view, the critical points of the Willmore functional, called {\it Willmore surfaces}, are thus also critical points of the Dirichlet energy of the Gauss map with respect to the induced metric $g$. Variations of $W$ in a fixed conformal class gives rise to a more general class of surfaces called {\it (conformal) constrained Willmore surfaces}.\\

Let $\vec{\mathbb{I}}_0$ denote the trace-free part of the second fundamental form, namely
\bes
\vec{\mathbb{I}}_0\;:=\;\vec{\mathbb{I}}\,-\,g\otimes\bH\:.
\ees
The Euler-Lagrange equation ({\it Willmore equation}) associated with the functional $W$ reads
\bes
\Delta_\perp\bH\,+\,g^{ik}g^{jl}\big(\vec{\mathbb{I}}_0\big)_{ij}\big\langle \big(\vec{\mathbb{I}}_0\big)_{kl}\,,\bH\big\rangle_{\R^m}\;=\;\vec{0}\:,
\ees
where $\Delta_\perp$ is the negative covariant Laplacian for the connection in the normal bundle.\\
Varying the Willmore functional under infinitesimal, smooth, compactly supported, {\it conformal} variations yields (cf. \cite{BPP}) the constrained Willmore equation
\be\label{cwe}
\Delta_\perp\bH\,+\,g^{ik}g^{jl}\big(\vec{\mathbb{I}}_0\big)_{ij}\big\langle \big(\vec{\mathbb{I}}_0\big)_{kl}\,,\bH\big\rangle_{\R^m}\;=\;g^{ik}g^{jl}\big(\vec{\mathbb{I}}_0\big)_{ij}q_{kl}\:,
\ee
where $q$ is a transverse traceless symmetric 2-covariant tensor\footnote{i.e. $q$ is divergence-free, where $(\text{div}\,q)_i:=\nabla^j q_{ji}$, and $\nabla^j$ denotes the covariant derivative.} with respect to the induced metric $g$. As is easily seen, $q$ plays the role of a Lagrange multiplier.\\
In \cite{Ri2} (see also \cite{Sc}), it is shown that solutions of the constrained Willmore equation are critical points of the Willmore energy constrained to a fixed conformal class. This notion clearly generalizes that of a Willmore surface, obtained via all smooth compactly supported infinitesimal variations (setting $q\equiv0$ in (\ref{cwe})). In the paper \cite{BR1}, the constrained Willmore equation arises as the limit of Palais-Smale sequences for the Willmore functional. \\

Minimal surfaces are examples of Willmore surfaces ; parallel mean curvature surfaces\footnote{parallel mean curvature surfaces satisfy $\pro d\bH\equiv\vec{0}$. They generalize to higher codimension the notion of constant mean curvature surfaces defined in $\R^3$. See the appendix.} are examples of constrained Willmore surfaces\footnote{{\it non}-minimal parallel mean curvature surfaces are however {\it not} Willmore.}. Not only is the Willmore energy invariant under reparametrization of the domain, but more remarkably, it is invariant under conformal transformations of $\R^m\cup\{\infty\}$.
Hence, the image of a [constrained] Willmore immersion through a conformal transformation is again a [constrained] Willmore immersion. It is thus no surprise that the class of Willmore immersions [resp. constrained Willmore immersions] is considerably larger than that of immersions whose mean curvature is minimal [resp. parallel], which are {\it not} preserved through conformal diffeomorphism. \\

In this paper, we will study the local regularity properties of constrained Willmore immersions. We impose no restriction on the codimension, and we allow the presence of point-singularities, called {\it branch points}. The only significant restriction will be imposed on the {\it Lagrange multiplier function} appearing in the Euler-Lagrange equation obtained through varying the Willmore energy in a fixed conformal class. We will demand that this multiplier be locally integrable (Willmore immersions and parallel mean curvature immersions have that feature). One of our goals will be to develop asymptotic expansions of the geometric quantities of the problem and seek conditions ensuring the {\it removability} of the branch point.\\

A {\it branch point} is a point where the immersion $\vec{\Phi}$ degenerates in the sense that $d\bp$ vanishes at that point. We focus on (conformal) locally Lipschitz and $W^{2,2}$ immersions $\bp:\di\setminus\{0\}\rightarrow\R^m$ with a branch point at the origin $0$, and regular away from the origin. {\it A priori}, at a branch point, the mean curvature is singular. 
We shall use the words {\it branch point} and {\it singularity} interchangeably. This is of course an abuse of language, as the immersion $\bp$ is not singular at a branch point. In fact, both $\bp$ and $d\bp$ are well-defined there. It is the immersive nature of $\bp$ which degenerates at a branch point.\\
In the context of this paper, the word {\it removability} is to be understood with care. To say that a branch point is removable does not mean that it is the result of some ``parametric illusion". Rather, it means that the map $\bp$ is smooth through the branch point, although it continues to fail to be an immersion at that point.\\

In the related works \cite{BR2} and \cite{BR3}, in collaboration with T. Rivi\`ere, the author delved deeper into the analysis of sequences of Willmore surfaces with uniformly bounded energy and non-degenerating conformal type. Our present work naturally borrows a lot from the techniques developed originally in \cite{BR3}.

\subsection{Main Results}\label{secloc}

The Willmore equation, which we will recall below, is a fourth-order nonlinear system of strongly couple equations (reducing to one scalar equation in codimension 1) for the immersion. In codimension 1, techniques have been developed to study its properties (cf. \cite{KS1} and the references therein). In higher codimension, new techniques were originally developed by Rivi\`ere in the seminal paper \cite{Ri1}, whereby the author shows that in a suitable conformal reparametrization, the Willmore equation can be written in divergence form with respect to flat local coordinates.  \\

\noindent
The constrained Willmore equation is the Willmore equation supplemented with a term involving both geometric quantities and an extraneous Lagrange multiplier, which in general does not depend on the geometric data of the problem. Our first task in this paper will be to show that in a suitable conformal reparametrization, the constrained Willmore equation can also be written in divergence form with respect to flat local coordinates.

\subsubsection{Reformulation of the problem}\label{refprob}

As we allow for branch points, we assume that the point-singularity lies at the origin, and we localize the problem by considering a map $\bp:\di\rightarrow\R^{m\ge3}$, which is an immersion of $\di\setminus\{0\}$, and satisfying
\begin{itemize}
\item[(i)] $\:\:\:\bp\in C^0(\di)\cap C^\infty(\di\setminus\{0\})\:;$
\item[(ii)] $\:\:\:\mathcal{H}^2\big(\bp(\di)\big)\,<\,\infty\:;$
\item[(iii)] $\displaystyle{\:\:\:\int_{\di}|\vec{\mathbb{I}}|^2_g\,d\text{vol}_g\,<\,\infty}\:.$
\end{itemize}
By a procedure detailed in \cite{BR3} (see also the references therein), it is possible to construct a parametrization $f$ of the unit-disk such that $\bp\circ f$ is {\it conformal}. To do so, one first extends $\bp$ to all of $\C\setminus\{0\}$ while keeping a bounded image and the second fundamental form square-integrable. One then shifts so as to have $\bp(0)=\vec{0}$, and inverts about the origin so as to obtain a complete immersion with square-integrable second fundamental form. Calling upon a result of Huber \cite{Hu} (see also \cite{MS}), one deduces that the image of the immersion is conformally equivalent to $\C$. Inverting yet once more about the origin finally gives the desired conformal immersion\footnote{which degenerates at the origin in a particular way, see (\ref{immas}).}, which we shall abusively continue to denote $\bp$. It has the aforementioned properties (i)-(iii), and moreover,
\bes
\bp(0)\,=\,\vec{0}\qquad\text{and}\qquad\bp(\di)\subset B^m_R(0)\quad\text{for some $\:0<R<\infty$}\:.
\ees
Hence, $\bp\in W^{1,\infty}\cap W^{2,2}(\di\setminus\{0\})$. Away from the origin, we define the {\it Gauss map} $\bn$ via
\bes
\bn\;=\;\star\,\dfrac{\px\bp\wedge\py\bp}{|\px\bp\wedge\py\bp|}\:\:,
\ees
where $\{x_1,x_2\}$ are standard Cartesian coordinates on the unit-disk $\di$, and $\star$ is the Euclidean Hodge-star operator. The immersion $\bp$ is conformal, i.e.
\be\label{confcond}
|\px\bp|\;=\;\text{e}^{\la}\;=\;|\py\bp|\qquad\text{and}\qquad\px\bp\cdot\py\bp\;=\;0\:,
\ee
where $\la$ is the conformal parameter. 
An elementary computation shows that
\be\label{eraser}
d\text{vol}_g\;=\;\text{e}^{2\la}dx\qquad\text{and}\qquad|\nabla\bn|^2dx\;=\;|d\bn|_g^2\,d\text{vol}_g\;=\;|\vec{\mathbb{I}}|^2_g\,d\text{vol}_g\:.
\ee
Hence, by hypothesis, we see that $\bn\in W^{1,2}(\di\setminus\{0\})$. In dimension two, the 2-capacity of isolated points is null, so we actually have $\bn\in W^{1,2}(\di)$. Rescaling if necessary, we shall henceforth always assume that 
\be\label{acheumeuneu}
\int_{\di}|\nabla\bn|^2\,dx\;<\;\eps_0\:,
\ee
where the adjustable parameter $\eps_0\equiv\eps_0(m)$ is chosen to fit our various needs (in particular, we will need it to be ``small enough" in Proposition \ref{morreydecay}). \\

For the sake of the following paragraph, we consider a conformal immersion $\bp:\di\rightarrow\R^m$, which is smooth across the unit-disk. We introduce the local coordinates $(x_1,x_2)$ for the flat metric on the unit-disk $\,\di=\big\{x=(x_1,x_2)\in{\R}^2\ ;\ x_1^2+x_2^2<1\big\}$. The operators $\nabla=(\px,\py)$, $\nabla^{\perp}=(-\py,\px)$, $\text{div}=\nabla\cdot\,$, and $\Delta=\nabla\cdot\nabla$  will be understood in these coordinates.
The conformal parameter $\la$ is defined as in (\ref{confcond}). We set
\be\label{lesvec}
\bej\::=\:\text{e}^{-\la}\pj\bp\qquad\text{for}\quad j\,\in\,\{1,2\}\:.
\ee
As $\bp$ is conformal, $\{\bex(x),\bey(x)\}$ forms an orthonormal basis of the tangent space $T_{\bp(x)}\bp(\di)$. Owing to the topology of $\di$, there exists for almost every $x\in\di$ a positively oriented orthonormal basis $\{\bn_1,\ldots,\bn_{m-2}\}$ of the normal space $N_{\bp(x)}\bp(\di)$, such that 
$\{\bex,\bey,\bn_1,\ldots,\bn_{m-2}\}$ forms a basis of $T_{\bp(x)}\R^m$. From the Pl\"ucker embedding, realizing the Grassmannian $Gr_{m-2}(\R^m)$ as a submanifold of the projective space of the $(m-2)^\text{th}$ exterior power $\,\mathbb{P}\big(\bigwedge^{m-2}\R^m\big)$, we can represent the Gauss map as the $(m-2)$-vector $\,\bn=\bigwedge_{\alpha=1}^{m-2}\bn_\alpha.$
Via the Hodge operator $\star\,$, we identify vectors and $(m-1)$-vectors in ${\R}^m$, namely:
\bes
\star\,(\bn\wedge \bAe_1)\,=\,\bAe_2\quad,\qquad\star\,(\bn\wedge \bAe_2)\,=\,-\,\bAe_1\quad,\qquad\star\,(\bex\wedge\bey)\,=\,\bn\:.
\ees
In this notation, the second fundamental form $\vec{\mathbb{I}}$, which is a symmetric 2-form on $T_{\bp(x)}\bp(\di)$  into  $N_{\bp(x)}\bp(\di)$, is expressed as
\bes
\vec{\mathbb{I}}\;=\;\sum_{\alpha,i,j}\ \text{e}^{-2\la}\,h^\alpha_{ij}\ \bn_\alpha\,dx_i\otimes dx_j\;\equiv\;\sum_{\alpha,i,j}\ h^\alpha_{ij}\ \bn_\alpha\,(\bAe_i)^\ast\otimes(\bAe_j)^\ast\:,
\ees
where
\bes
h^\alpha_{ij}\;=\;-\,\text{e}^{-\la}\,\bei\cdot\pj\bn_\alpha\:.
\ees\\
The mean curvature vector is
\bes
\bH\:=\:\sum_{\alpha=1}^{m-2}\,H^\alpha\,\bn_\alpha\:=\:\;\frac{1}{2}\,\sum_{\alpha=1}^{m-2}\,\big(h^\alpha_{11}+h^\alpha_{22}\big)\, \bn_\alpha\:,
\ees
and the Weingarten vector is
\bes
\bH_0\:=\:\sum_{\alpha=1}^{m-2}\,H_0^\alpha\,\bn_\alpha\:=\:\;\frac{1}{2}\,\sum_{\alpha=1}^{m-2}\,\big(h^\alpha_{11}-h^\alpha_{22}-2ih^\alpha_{12}\big)\, \bn_\alpha\:,
\ees
\noindent
As the tensor $q$ appearing in (\ref{cwe}) is traceless, symmetric, and divergence-free with respect to the conformal metric (\ref{eraser}), it is easily verified that
\bes
\px q_{11}\;=\;-\,\py q_{12}\qquad\text{and}\qquad \py q_{11}\;=\;\px q_{12}\:,
\ees
so that the function $f:=q_{11}+iq_{12}$ is anti-holomorphic. \\
Introducing $\bH_0$ and $f$ into (\ref{cwe}) enables us to recast it in the form
\be\label{wil1}
\Delta_\perp\bH\,+\,2\,\Re\big((\bH\cdot\bH^*_0)\bH_0\big)\:=\:\text{e}^{-2\la}\Re({\bH_0}f)\:,
\ee
with
\bes
-\,\Delta_\perp\bH\::=\:\text{e}^{-2\la}\,\pi_{\bn}\,\text{div}\,\pi_{\bn}\nabla\bH\:,
\ees
and $\pro$ is the projection onto the normal space spanned by $\{\bn_\alpha\}_{\alpha=1}^{m-2}$. We use the upper star notation to indicate complex conjugation.\\
To ease the notation, it is convenient to introduce the complex coordinates $z:=x_1+ix_2$ and $\bar{z}:=x_1-ix_2$. Then $\pz:=\frac{1}{2}(\px-i\py)$ and $\pzb:=\frac{1}{2}(\px+i\py)$. Similarly, we let
\bes
\bez\;:=\;\dfrac{1}{2}\,(\bex-i\bey)\qquad\text{and}\qquad \bezb\;:=\;\dfrac{1}{2}\,(\bex+i\bey)\:.
\ees
This notation is particularly helpful to recast the mean curvature and Weingarten vectors in simple forms, namely
\be\label{weinga}
\pzb\big(\text{e}^{\la}\bez\big)\;=\;\dfrac{\text{e}^{2\la}}{2}\,\bH\qquad\text{and}\qquad\pz\big(\text{e}^{-\la}\bez\big)\;=\;\dfrac{1}{2}\,\bH_0\:.
\ee

\noindent
The Willmore equation (\ref{wil1}) is a fourth-order nonlinear equation (with respect to the immersion $\bp$). With respect to the coefficients $H^\alpha$ of the mean curvature vector, it is actually a strongly coupled nonlinear system whose study is particularly challenging. Fortunately, in a conformal parametrization, it is possible to recast the left-hand side of (\ref{wil1}) in an equivalent, yet analytically more suitable form \cite{Ri1}. Namely, there holds on one hand
\begin{eqnarray}\label{wildiv}
&&\hspace{-.5cm}\Delta_\perp\bH\,+\,2\,\Re\big((\bH\cdot\bH^*_0)\bH_0\big)
\:=\:\dfrac{\text{e}^{-2\la}}{2}\,\text{div}\Big(\nabla\bH\,-\,3\,\pro\nabla\bH\,+\,\star\,(\nabla^\perp\bn\wedge\bH)\Big)\nonumber\\[1ex]
&&\hspace{2cm}=\:\:2\,\text{e}^{-2\la}\,\Re\Big[\pz\Big(\pzb\bH\,-\,3\,\pro\pzb\bH\,+\,i\star\big(\pzb\bn\wedge\bH\big)\Big)\Big]\:.\qquad
\end{eqnarray}
On the other hand, because $f$ is anti-holomorphic, (\ref{weinga}) gives us\footnote{we exclude the origin from the domain of validity of this identity, since we will in time allow for $f$ to be singular at the origin.}
\be\label{pety}
\pz\big(\text{e}^{-\la}f\,\bez\big)\;=\;\dfrac{1}{2}\,{\bH_0}f\qquad\text{on}\:\:\:D^2\setminus\{0\}\:.
\ee
Combining the latter to (\ref{wildiv}) yields
\begin{eqnarray*}
&&\Delta_\perp\bH\,+\,2\,\Re\big((\bH\cdot\bH^*_0)\bH_0\big)\,-\,\text{e}^{-2\la}\,\Re(\bH_0f\big)\nonumber\\[1ex]
&&\hspace{1cm}=\:\:2\,\text{e}^{-2\la}\,\Re\Big[\pz\Big(\pzb\bH\,-\,3\,\pro\pzb\bH\,+\,i\star\big(\pzb\bn\wedge\bH\big)\,-\,\text{e}^{-\la}f\,\bez\Big)\Big]\:.\qquad
\end{eqnarray*}
The constrained Willmore equation may thus be recast in the form
\be\label{confdiv}
\Re\Big[\pz\Big(\pzb\bH\,-\,3\,\pro\pzb\bH\,+\,i\star\big(\pzb\bn\wedge\bH\big)\,-\,\text{e}^{-\la}f\,\bez\Big)\Big]\;=\;\vec{0}\:.
\ee
In the $\{x_1,x_2\}$-variables, it is equivalently written
\be\label{confdiv2}
\text{div}\Big[\nabla\bH\,-\,3\,\pro\nabla\bH\,+\,\star\,(\nabla^\perp\bn\wedge\bH)\,-\,\text{e}^{-2\la}M_f\nabla^\perp\bp\Big]\;=\;\vec{0}\:,
\ee
with
\be\label{defmf}
M_f\;:=\;\left(\begin{array}{rc}-\Im(f)&\Re(f)\\[.5ex] \Re(f)&\Im(f)\end{array}\right)\:.
\ee

\medskip
This remarkable reformulation in divergence form of the Willmore equation is the starting point of our analysis. If there is a branch point at the origin, or if $f$ is singular at the origin, equation (\ref{confdiv2}) holds only away from the origin, on $\di\setminus\{0\}$. In particular, we can define the constant $\,\vec{\beta}_0\in\R^m$, called {\it first residue}, by
\be\label{residudu}
\vec{\beta}_0\,:=\,\dfrac{1}{4\pi}\int_{\partial\di}\vec{\nu}\cdot\Big(\nabla\bH\,-\,3\,\pro\nabla\bH\,+\,\star\,(\nabla^\perp\bn\wedge\bH)\,-\,\text{e}^{-2\la}M_f\nabla^\perp\bp)\Big)\:,
\ee
where $\vec{\nu}$ denotes the unit outward normal vector to $\partial\di$. We will see in Corollary \ref{Th3} that the first residue appears in the local asymptotic expansion of the mean curvature vector around the singularity.

\subsubsection{First regularity results}

We first state a result describing the regularity of the Gauss map around the point-singularity at the origin.
\begin{Prop}\label{Th1}
Let $\bp\in C^\infty(\di\setminus\{0\})\cap (W^{2,2}\cap W^{1,\infty})(\di)$ be a conformal constrained Willmore immersion of the punctured disk into $\R^m$ with integrable multiplier function $f$, and whose Gauss map $\bn$ lies in $W^{1,2}(\di)$. Then $\nabla^2\bn\in L^{2,\infty}(\di)$, and thus in particular $\nabla\bn$ is an element of $BMO$. Furthermore, $\bn$ satisfies the pointwise estimate
\bes
|\nabla\bn(x)|\;\lesssim\;|x|^{-\epsilon}\qquad\forall\:\:\epsilon>0\:.
\ees
If the order of degeneracy of the immersion $\bp$ at the origin is at least two\footnote{Roughly speaking, if $\nabla\bp(0)=\vec{0}$. The notion of ``order of degeneracy" is made precise below.}, then $\nabla\bn$ belongs to $L^\infty(\di)$. 
\end{Prop}

\medskip

A conformal immersion of $\di\setminus\{0\}$ into $\R^m$ such that $\nabla\bp$ and the Gauss map $\bn$ both extend to maps in $W^{1,2}(\di)$ has a distinctive behavior near the point-singularity located at the origin. One shows (cf. \cite{MS}, and Lemma A.5 in \cite{Ri2}) that there exists a positive integer $\theta_0$ with
\be\label{immas}
|\bp(x)|\;\simeq\;|x|^{\theta_0}\qquad\text{and}\qquad|\nabla\bp(x)|\;\simeq\;|x|^{\theta_0-1}\qquad\text{near the origin}\:.
\ee
In addition, there holds
\be\label{deffu}
\la(x)\;:=\;\dfrac{1}{2}\log\Big(\frac{1}{2}\,\big|\nabla\bp(x)\big|^2\Big)\;=\;(\theta_0-1)\log|x|\,+\,u(x)\:,
\ee
where $u\in W^{2,1}(\di)$ ; and one has 
\be\label{lambdas}
\left\{\begin{array}{lcl}\nabla\la\,\in\,L^2(\di)&,&\text{when $\,\theta_0=1$}\\[1ex]
|\nabla\la(x)|\,\lesssim\,|x|^{-1}\,\in\,L^{2,\infty}(\di)&,&\text{when $\,\theta_0\ge2\:.$}\end{array}\right.
\ee
The function $\,\text{e}^{-u(x)}\equiv|x|^{\theta_0-1}\text{e}^{-\la(x)}$ is continuous and strictly positive in a small neighborhood of the origin. \\
The integer $\theta_0$ is the density of the current $\bp_*[\di]$ at the image point $0\in\R^m$.\\

\noindent
When such a conformal immersion is Willmore on $\di\setminus\{0\}$, it is possible to refine the asymptotics (\ref{immas}). The following result describes the behavior of the immersion $\bp$ locally around the singularity at the origin. 

\begin{Prop}\label{Th2}
Let $\bp$ be as in Proposition \ref{Th1} with conformal parameter $\la$, and let $\theta_0$ be as in (\ref{immas}). 
There exists a constant vector $\bA=\bA^{1}+i\bA^{2}\in\R^2\otimes\R^m$ such that
\bes
\bA^{1}\cdot\bA^{2}\,=\,0\:\:,\qquad|\bA^{1}|\,=\,|\bA^{2}|\,=\,\theta_0^{-1}\lim_{x\rightarrow0}\,\dfrac{{e}^{\la(x)}}{|x|^{\theta_0-1}}\:\:,\qquad\pi_{\bn(0)}\bA\,=\,\vec{0}\:\:,
\ees
and
\be\label{rept}
\bp(x)\;=\;\Re\big(\bA\,{x}^{\theta_0-1}\big)+\vec{\xi}(x)\:,
\ee
with
\bes
\vec{\xi}(x)\,=\,\text{O}(|x|^{\theta_0-\epsilon})\:\:,\qquad\nabla\vec{\xi}(x)\,=\,\text{O}(|x|^{\theta_0-1-\epsilon})\:\:,\qquad\forall\:\:\epsilon>0\:.
\ees
\end{Prop}

\medskip
\noindent
The plane $span\{\bA^1,\bA^2\}$ is tangent to the surface at the origin. If $\theta_0=1$, this plane is actually $T_0\Sigma$. One can indeed show that the tangent unit vectors $\bej(0)$ spanning $T_0\Sigma$ (defined in (\ref{lesvec})) satisfy $\bej(0)=\bA^j/|\bA^j|$. In contrast, when $\theta_0\ge2$, the tangent plane $T_0\Sigma$ does not exist in the classical sense, and the vectors $\bej(x)$ ``spin" as $x$ approaches the origin. More precisely, $T_0\Sigma$ is the plane $span\{\bA^1,\bA^2\}$ covered $\theta_0$ times.

\subsubsection{Local asymptotic expansions}

Until now, no specific reference to the order of the anti-holomorphic multiplier function $f$ has been needed (other than the assumption that $f$ be integrable around the singularity at the origin). Our next result requires that the local behavior of $f$ at the origin be specified. To this end, we write
\be\label{fufu}
f\;=\;a_{\mu}\overline{z}^{\,\mu}+f_0\qquad\text{with}\:\:\quad\mu\ge-1\:\:,\:\:\:\:a_\mu\in\C\setminus\{0\}\:\:,\:\:\:\:f_0\in C^\infty(\di)\:.
\ee
We then combine the first residue $\vec{\beta}_0$ defined in (\ref{residudu}) and the vector $\bA$ from Proposition \ref{Th2} to define the {\it modified residue} 
\be\label{modudu}
\vec{\gamma}_0\;:=\;\vec{\beta}_0\,+\,\dfrac{1}{2}\,\delta_{\mu,\theta_0-2}\,\theta_0\,\text{e}^{-2u(0)}\Re\big(a_{\mu}\bA\big)\:,
\ee
where $u$ is the function in (\ref{deffu}).

\smallskip

\noindent
The next proposition describes the asymptotic behavior of the mean curvature vector near the singularity at the origin in terms of the modified residue. 

\begin{Prop}\label{Th3}
Let $\bp$ be as in Proposition \ref{Th1}, $\theta_0$ be as in (\ref{immas}), $\mu$ be as in (\ref{fufu}), and let $\vec{\gamma}_0$ be as in (\ref{modudu}). 
There holds locally around the origin
\be
\bH+\vec{\gamma}_0\log|z|\;=\;\Re\big(\vec{E}\big)+\text{O}(|z|^{1-a-\epsilon})\qquad\forall\:\epsilon>0\:,
\ee
where $\bE$ is a meromorphic function with a pole at the origin of order
\bes
\al\in\big\{\max\{0\,,\theta_0-\mu-2\},\ldots,\theta_0-1\big\}\:.
\ees
\end{Prop}

We may view the function $\bE$ from the previous proposition as a string of $m$ complex-valued functions $\{E_j\}_{j=1,\ldots,m}$, each of which is meromorphic and possibly has a pole at the origin of order at most $a$. This prompts us to defining the following decisive quantity. 

\begin{Dfi}\label{defresidu}
The {\it second residue} associated with the immersion $\bp$ at the origin is the $\mathbb{N}^{m}$-valued vector
\be\label{residudu2}
\vec{\gamma}\;=\;(\gamma_1, \ldots, \gamma_m)\qquad\text{with}\qquad \gamma_j\;:=\;-\,\dfrac{1}{2i\pi}\int_{\partial\di} d\log{E_j}\,\in\,\mathbb{N} \:.
\ee
\end{Dfi}
The importance of $\vec{\gamma}$ cannot be overstated: it controls the leading-order singular behavior of the mean curvature at the origin, as the following statement shows.

\begin{Th}\label{Th4}
Let $\bp$, $\theta_0$, $\mu$, $\vec{\gamma}_0$ be as Proposition \ref{Th3}, and let
$\vec{\gamma}$ be as in (\ref{residudu2}). Define
\bes
\al\;:=\;\max_{1\le j\le m}\gamma_j\,\in\big\{\max\{0\,,\theta_0-\mu-2\},\ldots,\theta_0-1\big\}\:.
\ees
The Gauss map satisfies
\bes
\nabla^{\theta_0+1-\al}\bn\,\in\,L^{2,\infty}(\di)\qquad\text{and thus}\qquad \nabla^{\theta_0-\al}\bn\,\in\,BMO(\di)\:.
\ees
Locally around the origin, the immersion has the asymptotic expansion
\bes
\bp=\Re\bigg(\bA\,z^{\theta_0}+\sum_{j=1}^{\theta_0-\al}\!\bB_j\,z^{\theta_0+j}+\bC_{\theta_0-\al}|z|^{2\theta_0}z^{-a}\bigg)-\bC|z|^{2\theta_0}\big(\log|z|^{\theta_0}-1\big)+\vec{\xi}\:,
\ees
where $\bB_j\in\C^m$ are constant vectors, while $\bA$ is as in Proposition \ref{Th2}. The constants $\bC_{\theta_0-a}$ and $\bC$ are
\bes
\bC_{\theta_0-a}\;:=\;\dfrac{\text{e}^{2u(0)}}{2\theta_0(\theta_0-a)}\,\vec{E}_a\qquad\text{and}\qquad \bC\;:=\;\dfrac{\text{e}^{2u(0)}}{2\theta_0^3}\,\vec{\gamma}_0\:,
\ees
where $\bE_a\in\C^m$ is a constant vector, and $\vec{\gamma}_0$ is the modified residue defined in (\ref{modudu}). The function $u$ is as in (\ref{deffu}).
The remainder $\vec{\xi}$ satisfies
\bes
\nabla^j\vec{\xi}\;=\;\text{O}(|z|^{2\theta_0-a+1-j-\epsilon})\qquad\forall\:j\in\{0,\ldots,\theta_0-a+1\}\:,\:\:\forall\:\,\epsilon>0\:.
\ees
The local behavior of the mean curvature follows accordingly: 
\bes
\bH\;=\;\Re\big(\vec{E}_az^{-a}\big)-\vec{\gamma}_0\log|z|+\,\vec{\eta}\:,
\ees
with
\bes
|z|^{\al-1+j}\nabla^j\vec{\eta}\;\in\,\bigcap_{p<\infty}L^p\qquad\forall\:j\in\{0,\ldots,\theta_0-\al\}\:,\:\:\forall\:\,\epsilon>0\:.
\ees
\end{Th}

\smallskip

The statement of Theorem \ref{Th4} reveals that in general, there is no reason to hope that the immersion should be smooth even when the multiplier function is. In fact, there holds

\bes
\bp\,\in\bigcap_{p<\infty}\left\{\begin{array}{lcl}
W^{2,p}&\:,\:&\theta_0=1\\[1ex]
W^{\theta_0+2-a,p}&\:,\:&\theta_0\ge2\:.
\end{array}\right.
\ees
In the worst case scenario, $a=\theta_0-1$. Then the immersion might be as little regular as $C^{1,\alpha}(\di)$ (when $\theta_0=1$). The example of the inverted catenoid (cf. \cite{BR3}), which is a Willmore surface, hence with $f\equiv0$, reveals that the surface may not be $C^{1,1}$. When the singularity has a higher order, $\theta_0\ge2$, the immersion is always $C^{2,\alpha}$ for all $\alpha<1$. 

\subsubsection{Vanishing residues: improved regularity}

At the other end of the spectrum, the best case scenario occurs when $a=0$. This is however not always possible to achieve: it is necessary that the multiplier function $f$ decays fast enough at the origin, namely, that $\mu\ge\theta_0-2$. The next result details what happens when both the modified residue and the second residue vanish (which is tantamount to $a=0$). 

\begin{Th}\label{Th44}
Under the hypotheses of Theorem \ref{Th4}, if the modified residue $\vec{\gamma}_0$ and if the second residue $\vec{\gamma}$ both vanish\footnote{if $\theta_0=1$, the second residue $\vec{\gamma}$ automatically vanishes.}, then there holds
\begin{itemize}
\item[(i)] if $\theta_0<\mu+2$, the immersion is smooth throughout the branch point. 
\item[(ii)] if $\theta_0=\mu+2$, the mean curvature vector satisfies 
\bes
\bH\,\in\,W^{2,(2,\infty)}(\di)\,\subset\bigcap_{\alpha\in[0,1)}C^{0,\alpha}(\di)\:.
\ees
Furthermore, there holds
\bes
\bp\,\in\left\{\begin{array}{l}
W^{4,(2,\infty)}\\[1ex]
\bigcap_{p<\infty}W^{\theta_0+2,p}
\end{array}\right.\quad\text{and}\qquad
\bn\,\in\left\{\begin{array}{lcl}
W^{3,(2,\infty)}&,&\theta_0=1\\[1ex]
\bigcap_{p<\infty}W^{\theta_0+1,p}&,&\theta_0\ge2\:.
\end{array}\right.
\ees
In particular,
\bes
\bp\,\in\,C^{\theta_0+1,\alpha}(\di)\qquad\text{and}\qquad \bn\,\in C^{\theta_0,\alpha}(\di)\qquad\forall\:\alpha\in[0,1)\:.
\ees
\end{itemize}
\end{Th}

In the special case when the origin is a regular point (i.e. when the constrained Willmore equation (\ref{wil1}) holds on the whole unit disk, the regularity of the solution follows from Theorem \ref{Th44} by setting $\theta_0=1$\footnote{indeed, immersions which have either a regular point at the origin or a branch point of order one behave analogously, namely $|\nabla\bp|$ is bounded from above and below at the origin.}. One readily checks that both residues $\vec{\gamma}_0$ and $\vec{\gamma}$ vanish at a regular point. Whence we deduce

\begin{Co}\label{Co1}
When the origin is a regular point (i.e. when constrained Willmore equation (\ref{wil1}) holds on the whole unit disk), we have
\begin{itemize}
\item[(i)] if the multiplier function $f$ is regular at the origin (i.e. if $\mu\ge0$), then the immersion is smooth across the singularity.
\item[(ii)] if the multiplier function $f$ is singular at the origin (i.e. if $\mu=-1$), then the immersion lies in $C^{2,\alpha}(\di)$ for all $\alpha\in[0,1)$. 
\end{itemize}
\end{Co}

\subsubsection{Surfaces of specific types}

Naturally, the statement of Theorem \ref{Th44} applies to constrained Willmore surfaces of specific types. When $f\equiv0$, the immersion is called Willmore. The behavior of a Willmore surface near a point-singularity was settled in \cite{BR3}, whose results are recovered by setting $f\equiv0$ in the statement of Theorem \ref{Th4}. Furthermore, when both the first residue $\vec{\beta}_0$ and the second residue $\vec{\gamma}$ vanish, the Willmore immersion is smooth across the singularity at the origin. \\

The second important subclass of constrained Willmore surfaces are the parallel mean curvature surfaces, whose multiplier function is constrained by the geometry of the problem, namely $f=2\text{e}^{2\la}\bH\cdot\bH_0^*$. It is not difficult to verify that a parallel mean curvature surface with a branch point at the origin has always vanishing residues, and moreover that $\mu\ge\theta_0-1$. Theorem \ref{Th44} then guarantees that the immersion is smooth across the point-singularity.\\

 For clarity, the above results are summarized in the following corollary. 

\begin{Co}\label{Co2}
Let $\bp$ satisfy the hypotheses of Proposition \ref{Th1}.
\begin{itemize}
\item[(i)] If $\bp$ is Willmore (i.e. $f\equiv0$), then $\bp$ is smooth across regular points. It is smooth across a point-singularity provided the residues $\vec{\beta}_0$ and $\vec{\gamma}$ vanish.  
\item[(ii)] If $\bp$ has parallel mean curvature (i.e. $f=2\text{e}^{2\la}\bH\cdot\bH_0^*$), then $\bp$ is smooth across regular points and across branch points alike. 
\end{itemize}
\end{Co}


\subsubsection{Open questions}

The first open question posed by our results is that of sharpness. In the case of Willmore immersions (when the multiplier function $f$ identically vanishes), it was firmly established in \cite{BR3} that the results are sharp. This is much less clear when the multiplier function is not zero. The difficulty essentially stems from the lack of relevant examples.\\ 


The second, and most analytically challenging, open question deals with studying the local behavior of a constrained Willmore immersion with multiplier function degenerating at the origin with order strictly greater than one (i.e. when $f$ is not integrable). Even when the origin is assumed to be a geometrically regular point (that is, when the origin is not a branch point), it is already unclear how to proceed. The author conjectures that if the multiplier function has a pole of order strictly greater than 3 anywhere, then the associated constrained Willmore equation can have no solution $\bp\in W^{2,2}$. The case of a multiplier function having a pole of order exactly 2 seems to be the most problematic\footnote{although the author knows no example displaying this feature, there seems to be no analytical reason why such examples should be ruled out.}. The heart of the matter lies within the conformal conservative Willmore system (\ref{sysSR}) which cannot be extended through the origin when the multiplier function is not integrable. It is thus not possible to reduce the problem into a subcritical form. For the same reason, if one assumes that the immersion lies in the Sobolev-Lorentz\footnote{cf. the Appendix of \cite{BR1} for further details on Sobolev-Lorentz spaces.} space $W^{2,(2,1)}(\di)$, which is slightly smaller than $W^{2,2}(\di)$, it follows that the Gauss map $\bn$ has a well-defined limit at the origin, and mimicking the analysis presented here then becomes possible (this essentially amounts to assuming {\it de facto} that the problem is subcritical).

\reset

\section{Proof of the Theorems}

\subsection{Fundamental results and reformulation}

We place ourselves in the situation described in subsection \ref{refprob} of the Introduction. Namely, we have a constrained Willmore immersion $\bp$ on the punctured disk which degenerates at the origin in such a way that
\bes
|\bp(x)|\;\simeq\;|x|^{\theta_0}\qquad\text{and}\qquad|\nabla\bp(x)|\;=\;\sqrt{2}\,\text{e}^{\la(x)}\;\simeq\;|x|^{\theta_0-1}\:,
\ees
for some $\theta_0\in\N\setminus\{0\}$.\\[1ex]
The multiplier function $f$ is assumed to be integrable on the unit-disk. As it is anti-holomorphic, there holds $f=\text{O}(|x|^{-1})$ about the origin.

\medskip

Amongst the analytical tools available to the study of constrained Willmore immersions with square-integrable second fundamental form, an important one is certainly the so-called {\it $\eps$-regularity estimate}. The version appearing in Theorem 2.10 and Remark 2.11 from \cite{KS2} states that there exists $\eps_0>0$ such that, if
\bes
\int_{D^2}|\nabla\bn|^2\,dx\;<\;\eps_0\:,
\ees
then there holds
\begin{eqnarray}\label{epsreg}
\hspace{-.8cm}\Vert\text{e}^{-\la}\nabla\bn\Vert^2_{L^\infty(B^g_\sigma)}&&\\[0ex]
&&\hspace{-1.5cm}\;\le\;C_0\bigg[\Vert\text{e}^{-\la}\Re(\bH_0f)\Vert_{L^2(B^g_{2\sigma})}+\dfrac{1}{\sigma^2}\,\Vert\nabla\bn\Vert_{L^2(B^g_{2\sigma})}\bigg]\Vert\nabla\bn\Vert_{L^2(B^g_{2\sigma})}\:,\nonumber
\end{eqnarray}
where $B_\sigma^g$ is any geodesic disk of radius $\sigma$ for the induced metric $g=\bp^*g_{\R^m}$ with $B^g_{2\sigma}\subset\Om:=D^2\setminus\{0\}$ ; and $C_0$ is a universal constant. As always, $\la$ denotes the conformal parameter.\\[1.5ex]
The $\eps$-regularity enables us to obtain the following result, decisive to the remainder of the argument.

\begin{Lma}\label{delta}
The function $\,\delta(r):=r\sup_{|x|=r}|\nabla\bn(x)|\,$ satisfies
\bes
\lim_{r\searrow0}\:\delta(r)\,=\,0\hspace{1cm}\text{and}\hspace{1cm}\int_{0}^{1}\delta^2(r)\,\dfrac{dr}{r}\,<\,\infty\:.
\ees
\end{Lma}
$\textbf{Proof.}$ 
From (\ref{eraser}) and (\ref{immas}), the metric $g$ satisfies 
\bes
g_{ij}(x)\;\simeq\;|x|^{2(\theta_0-1)}\delta_{ij}\qquad\text{on}\:\:\:\:D_{2r}(0)\setminus D_{r/2}(0)\qquad\forall\:\,r\in(0\,,1/2)\:.
\ees
A simple computation then shows that
\be\label{barbitruc}
B^g_{2cr^{\theta_0}}(x)\;\subset\;D_{2r}(0)\setminus D_{r/2}(0)\qquad\forall\:\,x\in \partial D_{r}(0)\:,
\ee
where $\,0<2\theta_0\,c<1-2^{-\theta_0}$.\\
Since the metric $g$ does not degenerate away from the origin, given $0<r<1/2$, we can always cover the flat circle $\partial D_{r}(0)$ with finitely many metric disks:
\bes
\partial D_r(0)\;\subset\;\bigcup_{j=1}^{N} B^g_{cr^{\theta_0}}(x_j)\qquad\text{with}\quad x_j\in\partial D_r(0)\:.
\ees
Hence, per the latter, (\ref{epsreg}), and (\ref{barbitruc}), we obtain that for some $y\in\partial D_r(0)$ there holds
\begin{eqnarray}\label{barbibulle}
r^2\sup_{|x|=r}|\nabla\bn(x)|^2&\simeq&r^{2\theta_0}\!\sup_{|x|=r}\,\big|\text{e}^{-\la(x)}\nabla\bn(x)\big|^2\;\;\leq\;\;r^{2\theta_0}\Vert\text{e}^{-\la}\nabla\bn\Vert^2_{L^\infty(B^g_{cr^{\theta_0}}(y))}\nonumber\\[0ex]
&&\hspace{-3.5cm}\lesssim\:\: r^{2\theta_0}\bigg[\Vert\text{e}^{-\la}\Re(\bH_0f)\Vert_{L^2(B^g_{2cr^{\theta_0}}(y))}+\dfrac{1}{r^{2\theta_0}}\,\Vert\nabla\bn\Vert_{L^2(B^g_{2cr^{\theta_0}}(y))}\bigg]\Vert\nabla\bn\Vert_{L^2(B^g_{2cr^{\theta_0}}(y))}\nonumber\\[1ex]
&&\hspace{-3.5cm}\lesssim\:\: r^{2\theta_0}\bigg[\Vert\text{e}^{-2\lambda}f\Vert_{L^\infty(B^g_{2cr^{\theta_0}}(y))}+\dfrac{1}{r^{2\theta_0}}\bigg]\Vert\nabla\bn\Vert^2_{L^2(B^g_{2cr^{\theta_0}}(y))}\nonumber\\[1ex]
&&\hspace{-3.5cm}\lesssim\:\: r^{2\theta_0}\bigg[\Vert\text{e}^{-2\lambda}f\Vert_{L^\infty(D_{2r}(0)\setminus D_{r/2}(0))}+\dfrac{1}{r^{2\theta_0}}\bigg]\Vert\nabla\bn\Vert^2_{L^2(D_{2r}(0)\setminus D_{r/2}(0))}\nonumber\\[1ex]
&&\hspace{-3.5cm}\lesssim\:\:\big(1+r\big)\,\Vert\nabla\bn\Vert^2_{L^2(D_{2r}(0)\setminus D_{r/2}(0))}\:.
\end{eqnarray}
We have used the facts that near the origin $\text{e}^{\la}\simeq|x|^{\theta_0-1}$, that $|f|\lesssim|x|^{-1}$, and that $\text{e}^{\la}|\bH_0|\lesssim|\nabla\bn|$. Whence,
\bes
\delta(r)\;:=\;r\sup_{|x|=r}|\nabla\bn(x)|\;\lesssim\;\Vert\nabla\bn\Vert_{L^2(D_{2r}(0)\setminus D_{r/2}(0))}\:.
\ees
As $\nabla\bn$ is square-integrable by hypothesis, letting $r$ tend to zero in the latter yields the first assertion. \\
The second assertion follows from (\ref{barbibulle}), namely, 
\bes
\int_{0}^{1/2}\delta^2(r)\,\dfrac{dr}{r}\;\lesssim\;\int_{0}^{1/2}\Vert\nabla\bn\Vert^2_{L^2(D_{2r}(0)\setminus D_{r/2}(0))}\,\dfrac{dr}{r}\;=\;\log(4)\;\Vert\nabla\bn\Vert^2_{L^2(\di)}\:,
\ees
which is by hypothesis finite.\\[-3ex]

$\hfill\blacksquare$ \\

We obtain from Lemma \ref{delta} that
\be\label{es2}
r^{\theta_0}\!\sup_{|x|=r}|\bH(x)|\:\lesssim\:\delta(r)\qquad\text{and}\qquad r^{\theta_0}\!\sup_{|x|=r}|\bH_0(x)|\:\lesssim\:\delta(r)\:.
\ee
As $\bn\wedge\bH=\vec{0}$, the constrained Willmore equation (\ref{confdiv2}) may be alternatively written
\be\label{daka0}
\text{div}\Big(-2\,\nabla\bH\,+\,3\,\pi_T{\nabla\bH}\,-\,\star\,(\bn\wedge\pi_T\nabla^\perp\bH)-\,\text{e}^{-2\la}M_f\nabla^\perp\bp\Big)\;=\;\vec{0}\:,
\ee
where $\pi_T:=\text{id}-\pro$ denotes projection onto the tangent space. 
Using the fact that $\bH$ is normal, a simple computation reveals that
\be\label{daka1}
\pi_T\nabla\bH\;=\;-\sum_{j=1,2}\big(\bH\cdot\pro\nabla\vec{e}_{j}\big)\,\vec{e}_{j}\;=\;-\,|\bH|^2\nabla\bp\,-\big(\bH\cdot M_{\bH_0}\big)\nabla^\perp\bp\:,
\ee
with
\bes
M_{\bH_0}\;:=\;\left(\begin{array}{cr}\Im(\bH_0)&\Re(\bH_0)\\[.5ex] \Re(\bH_0)&-\Im(\bH_0)  \end{array}\right)\:.
\ees
From this and the elementary identities
\bes
\star\,(\bn\wedge\nabla\bp)\;=\;-\,\nabla^\perp\bp\qquad\text{and}\qquad \star(\bn\wedge\nabla^\perp\bp)\;=\;\nabla\bp\:,
\ees
we obtain
\be\label{daka2}
\star\,(\bn\wedge\pi_T\nabla^\perp\bH)\;=\;-\,|\bH|^2\nabla\bp\,+\big(\bH\cdot M_{\bH_0}\big)\nabla^\perp\bp\:.
\ee
Combining (\ref{daka1}) and (\ref{daka2}) into (\ref{daka0}) yields
\be\label{newill}
\text{div}\bigg(\nabla\bH\,+\,|\bH|^2\nabla\bp\,+\,2\big(\bH\cdot M_{\bH_0}\big)\nabla^\perp\bp\,+\,\dfrac{1}{2}\,\text{e}^{-2\la}M_f\nabla^\perp\bp\bigg)=\;\vec{0}\:.
\ee
Observe that owing to (\ref{es2}), there holds
\bes
|x|^{\theta_0}\Big(|\bH|^2\nabla\bp\,+\,2\big(\bH\cdot M_{\bH_0}\big)\nabla^\perp\bp\Big)\;\lesssim\;\dfrac{\delta^2(|x|)}{|x|}\:.
\ees
On the other hand, 
\bes
|x|^{\theta_0}\text{e}^{-2\la}M_f\nabla^\perp\bp\;\lesssim\;|x||f|\;<\;\infty\:.
\ees
Since $|x|^{-1}\delta(|x|)$ and $|x|^{\theta_0-1}\bH$ both lie in $L^2$, we may apply Lemma \ref{GW} to (\ref{newill}) with $a=\theta_0$ so as to obtain that
\be\label{es3}
|x|^{\theta_0}\nabla\bH\,\in\,L^2(\di)\:.
\ee
This observation shall be helpful in the sequel. \\

Equation (\ref{wildiv}) implies that for any disk $D_\rho(0)$ of radius $\rho$ centered on the origin and contained in $\Om:=\di\setminus\{0\}$, there holds for all $\rho\in(0,1)$: 
\be\label{residue}
\int_{\partial D_\rho(0)} \vec{\nu}\cdot\Big(\nabla\bH-3\,\pro\nabla\bH+\star\,(\nabla^\perp\bn\wedge\bH)-\text{e}^{-2\la}M_f\nabla^\perp\bp\Big)\,=\,4\pi\vec{\beta}_0\:,
\ee
where $\vec{\beta}_0$ is the residue defined in (\ref{residudu}). 
An elementary computation shows that
\bes
\int_{\partial D_{\rho}(0)} \vec{\nu}\cdot\nabla\log|x|\:=\:2\pi\qquad\forall\:\rho>0\:.
\ees
Thus, upon setting
\be\label{a1}
\bX\;:=\;\nabla\bH\,-\,3\,\pro\nabla\bH\,+\,\star\,(\nabla^\perp\bn\wedge\bH)\,-\text{e}^{-2\la}M_f\nabla^\perp\bp\,-\,2\vec{\beta}_0\nabla\log|x|\:,
\ee
we find
\bes
\text{div}\,\bX=0\quad\:\:\text{on}\:\:\:\Om=\di\setminus\{0\}\qquad\text{and}\hspace{.5cm}\int_{\partial D_\rho(0)}\vec{\nu}\cdot\bX\,=\,{0}\:\:\quad\forall\:\rho\in(0,1)\:.
\ees
As $\bX$ is smooth away from the origin, the Poincar\'e Lemma implies now the existence of an element $\bL\in C^{\infty}(\Om)$, defined up to an additive constant, such that
\be\label{a2}
\bX\;=\;\nabla^\perp\bL\qquad\text{in}\:\:\:\Om\:.
\ee
We deduce from Lemma \ref{delta} and (\ref{es2})-(\ref{a2}) that
\be\label{xx1}
|x|^{\theta_0}\nabla\bL\,\in\,L^2(\di)\:.
\ee
A classical Hardy-Sobolev inequality gives the estimate
\be\label{estimL}
\theta_0^2\int_{\di}|x|^{2(\theta_0-1)}|\bL|^2\,dx\:\le\:\int_{\di}|x|^{2\theta_0}|\nabla\bL|^2\,dx\;+\;\theta_0\!\int_{\partial \di}|\bL|^2\:, 
\ee
which is a finite quantity, owing to (\ref{xx1}) and to the smoothness of $\bL$ away from the origin. The immersion $\bP$ has near the origin the asymptotic behavior $\,|\nabla\bp(x)|\simeq|x|^{\theta_0-1}$. Hence (\ref{estimL}) yields that 
\be\label{xx2}
\bL\cdot\nabla\bP\:,\;\bL\wedge\nabla\bP\;\in\;L^2(\di)\:.
\ee

We next set $\,\vec{\Gamma}(x):=2\vec{\beta}_0\log|x|$. 
In the paper \cite{BR1} (cf. Lemma A.2), two identities are derived:
\bes
\left\{\begin{array}{lcl}
\nabla\bP\cdot\big(\nabla^\perp\bL\,+\nabla\vec{\Gamma}\,+\text{e}^{-2\la}M_f\nabla^\perp\bp\big)&=&0\\[1ex]
\nabla\bP\wedge\big(\nabla^\perp\bL\,+\nabla\vec{\Gamma}\,+\text{e}^{-2\la}M_f\nabla^\perp\bp\big)&=&-\,2\,\nabla\bP\wedge\nabla\bH\:.\end{array}\right.
\ees
The conformality of $\bp$ yields easily that
\bes
\nabla\bp\cdot\big(M_f\nabla^\perp\bp\big)\;=\;0\;=\;\nabla\bp\wedge\big(M_f\nabla^\perp\bp\big)\:.
\ees
Whence, we find
\be\label{constraints}
\left\{\begin{array}{lcl}
\nabla\bP\cdot\big(\nabla^\perp\bL\,+\nabla\vec{\Gamma}\big)&=&0\\[1ex]
\nabla\bP\wedge\big(\nabla^\perp\bL\,+\nabla\vec{\Gamma}\big)&=&-\,2\,\nabla\bP\wedge\nabla\bH\:.\end{array}\right.
\ee
We find useful to define the functions $g$ and $\bG$ via
\be\label{sys1}
\left\{\begin{array}{rclcrclcl}
\Delta g&=&\nabla\vec{\Gamma}\cdot\nabla\bp&,\quad&\Delta\bG&=&\nabla\vec{\Gamma}\wedge\nabla\bp&\quad&\text{in}\:\:\: \di\\[1ex]
g&=&0&,\quad&\bG&=&\vec{0}&\quad&\text{on}\:\:\:\partial \di\:.
\end{array}\right.
\ee
Since $\,|\nabla\bP(x)|\simeq|x|^{\theta_0-1}\,$ near the origin and $\vec{\Gamma}$ is the fundamental solution of the Laplacian, by applying Calderon-Zygmund estimates to (\ref{sys1}), we find\footnote{The weak-$L^2$ Marcinkiewicz space $L^{2,\infty}(\di)$ is defined as those functions $f$ which satisfy $\:\sup_{\alpha>0}\alpha^2\Big|\big\{x\in \di\,;\,|f(x)|\ge\alpha\big\}\Big|<\infty$. In dimension two, the prototype element of $L^{2,\infty}$ is $|x|^{-1}\,$. The space $L^{2,\infty}$ is also a Lorentz space, and in particular is a space of interpolation between Lebesgue spaces, which justifies the first inclusion in (\ref{regg}). See \cite{He} for details.} 
\be\label{regg}
\nabla^2g\;,\;\nabla^2\bG\:\;\in\:\left\{\begin{array}{lcl}L^{2,\infty}(\di)&,&\theta_0=1\\[1.5ex]
BMO(\di)&,&\theta_0\ge2\:.
\end{array}\right.
\ee
Accounted into (\ref{sys1}), the latter yield that there holds in $\Om$\,:
\be\label{div1}
\left\{\begin{array}{rcl}
\text{div}\big(\bL\cdot\nabla^\perp\bP\,-\,\nabla g\big)&=&0\\[1.5ex]
\text{div}\big(\bL\wedge\nabla^\perp\bP\,-\,2\,\bH\wedge\nabla\bP\,-\,\nabla\bG\big)&=&\vec{0}\:,\end{array}\right.
\ee
where we have used the fact that
\bes
\Delta\bP\wedge\bH\;=\;2\,\text{e}^{2\la}\bH\wedge\bH\;=\;\vec{0}\:.
\ees
Note that the terms under the divergence symbols in (\ref{div1}) both belong to $L^2(\di)$, owing to (\ref{xx2}) and (\ref{regg}). The distributional equations (\ref{div1}), which are {\it a priori} to be understood on $\Om=\di\setminus\{0\}$, may thus be extended to all of $\di$. Indeed, a classical result of Laurent Schwartz states that the only distributions supported on $\{0\}$ are linear combinations of derivatives of the Dirac delta mass. Yet, none of these (including delta itself) belongs to $W^{-1,2}$. We shall thus understand (\ref{div1}) on $\di$. 
It is not difficult to verify (cf. Corollary IX.5 in [DL]) that a divergence-free vector field in $L^2(\di)$ is the curl of an element in $W^{1,2}(\di)$. 
We apply this observation to (\ref{div1}) so as to infer the existence of two functions\footnote{$S$ is a scalar while $\bR$ is $\bigwedge^2(\R^m)$-valued.} $S$ and of $\bR$ in the space $W^{1,2}(\di)\cap C^\infty(\Om)$, with
\bes
\left\{\begin{array}{rclll}
\nabla^\perp S&=&\bL\cdot\nabla^\perp\bP\,-\,\nabla g&&\\[1ex]
\nabla^\perp\bR&=&\bL\wedge\nabla^\perp\bP\,-\,2\,\bH\wedge\nabla\bP\,-\,\nabla\bG\:.
\end{array}\right.
\ees
According to Lemma A.1 in \cite{BR3}, the functions $S$ and $\bR$ satisfy on $\di$ the following system of equations, called {\it conservative conformal Willmore system}\footnote{refer to the Appendix for the notation and the operators used.}:
\be\label{sysSR}
\left\{\begin{array}{rclll}
-\,\Delta S&=&\nabla(\star\,\bn)\cdot\nabla^\perp\bR\,+\,\text{div}\big((\star\,\bn)\cdot\nabla\bG\big)&&\\[1.5ex]
-\,\Delta\bR&=&\nabla(\star\,\bn)\bul\nabla^\perp\bR\;-\;\nabla(\star\,\bn)\cdot\nabla^\perp S\\[.75ex]
&&\hspace{2.15cm}\,+\:\,\text{div}\big((\star\,\bn)\bul\nabla\bG\,+\,\star\,\bn\,\nabla g\big)\:.
\end{array}\right.
\ee
Not only is this system independent of the codimension (which enters the equations in the guise of the operators $\star$ and $\bul$), but it further displays two fundamental advantages. Analytically, (\ref{sysSR}) is uniformly elliptic. This is in sharp contrast with the constrained Willmore equation (\ref{wil1}) whose leading order operator $\Delta_{\perp}$ degenerates at the origin, owing to the presence of the conformal factor $\text{e}^{\la(x)}\simeq|x|^{\theta_0-1}$\,. Structurally, the system (\ref{sysSR}) is in divergence form. We shall in the sequel capitalize on this remarkable feature to develop arguments of ``integration by compensation". 
\medskip

The conservative conformal Willmore system (\ref{sysSR}) is also to be supplemented with the following important identity, also proved in Lemma A.1 of \cite{BR3}, namely 
\be\label{delphi}
-\,2\,\Delta\bP\;=\;\big(\nabla S-\nabla^\perp g\big)\cdot\nabla^\perp\bP\,-\,\big(\nabla\bR-\nabla^\perp\bG\big)\bul\nabla^\perp\bP\:.
\ee

\subsection{Preliminary estimates: from critical to subcritical}

{\it A priori}\, since $\bn$, $S$, and $\bR$ are elements of $W^{1,2}$, the leading terms on the right-hand side of the conservative conformal Willmore system (\ref{sysSR}) are critical. This difficulty can nevertheless be bypassed using the fact that the $W^{1,2}$-norm of the Gauss map $\bn$ is chosen to be small enough (cf. (\ref{acheumeuneu})). More precisely, in \cite{BR3} (cf. Proposition A.1), the following result is established  

\begin{Prop}\label{morreydecay}
Let $u\in W^{1,2}(\di)\cap C^2(\di\setminus\{0\})$ satisfy the equation
\bes
-\,\Delta u\:=\:\nabla b\cdot\nabla^\perp u+\text{div}\,(b\,\nabla h)\qquad\quad\text{on}\:\:\:\:\di\:,
\ees
where $\,h\in W_0^{2,(2,\infty)}(\di)$, and moreover
\bes
b\,\in\,W^{1,2}\cap L^{\infty}(\di)\qquad\text{with}\qquad\Vert\nabla b\Vert_{L^2(\di)}\;<\;\eps_0\:,
\ees
for some $\eps_0$ chosen to be ``small enough". Then
\bes
\nabla u\,\in\,L^p(\di)\qquad\text{for some $\,p>2$}\:.
\ees
\end{Prop}

\medskip
\noindent
Owing to (\ref{regg}) and to (\ref{acheumeuneu}), this result applies in particular to the system (\ref{sysSR}) and yields
\be\label{gradSinLp}
\nabla S\:,\:\nabla\bR\:\in\:L^{p}(\di)\qquad\text{for some}\:\:p>2\:.
\ee
Hence, as $\,|\nabla\bP(x)|\simeq\text{e}^{\la(x)}\simeq|x|^{\theta_0-1}$ around the origin, using (\ref{regg}) and (\ref{gradSinLp}), we may apply Proposition \ref{CZpondere} with the weight $|\mu|=\text{e}^{\la}$ and $a=\theta_0-1$ to the equation (\ref{delphi}) so as to conclude that
\bes
(\px+i\,\py)\bp(x)\;=\;\bPe(\overline{x})+\,\text{e}^{\la(x)}\bT(x)\:,
\ees
where $\bPe$ is a $\C^m$-valued polynomial of degree at most $(\theta_0-1)$, and $\vec{T}(x)=\text{O}\big(|x|^{1-\frac{2}{p}-\epsilon}\big)$ for every $\epsilon>0$. Because $\,\text{e}^{-\la}\nabla\bp\,$ is a bounded function, we deduce more precisely that $\bPe(\overline{x})=\theta_0\bA^*\,\overline{x}^{\,\theta_0-1}$, for some constant vector $\bA\in\C^m$\, (we denote its complex conjugate by $\bA^*$), so that
\bes
\nabla\bp(x)\;=\;\left(\begin{array}{c}\Re\\[.5ex]-\,\Im\end{array}\right)\big(\theta_0\,\vec{A}\;{x}^{\theta_0-1}\big)\,+\,\text{e}^{\la(x)}\bT(x)\:.
\ees
Equivalently, switching to the complex notation, there holds
\be\label{locexphi}
\pz\bp\;=\;\dfrac{\theta_0}{2}\bA\,z^{\theta_0-1}+\,\text{O}\big(|z|^{\theta_0-\frac{2}{p}-\epsilon}\big)\quad\forall\:\epsilon>0\:.
\ee
We write $\bA=\bA^1+i\bA^2\in\R^2\otimes\R^m$. 
The conformality condition on $\bp$ shows easily that $\bA\cdot\bA=0$, whence
\bes
|\bA^{1}|\;=\;|\bA^{2}|\qquad\text{and}\qquad\bA^{1}\cdot\bA^{2}\;=\;0\:.
\ees
Yet more precisely, as $|\nabla\bp|^2=2\,\text{e}^{2\la}$, we see that
\bes
|\bA^{1}|\;=\;|\bA^{2}|\;=\;\dfrac{1}{\theta_0}\,\lim_{z\rightarrow0}\,\dfrac{\text{e}^{\la(z,\bar{z})}}{|z|^{\theta_0-1}}\,\in\:\:]0\,,\infty[\:.
\ees
Because $\bp(0)=\vec{0}$, we obtain from (\ref{locexphi}) the local expansion
\bes
\bp(z,\bar{z})\;=\;\Re\big(\bA\,z^{\theta_0}\big)\,+\,\text{O}\big(|z|^{\theta_0+1-\frac{2}{p}-\epsilon}\big)\:.
\ees
On the other hand, from $\pro\nabla\bp\equiv\vec{0}$, we deduce from (\ref{locexphi}) that
\bes
\pi_{\bn}\bA\;=\;\text{O}\big(|z|^{1-\frac{2}{p}-\epsilon}\big)\qquad\:\:\:\forall\:\epsilon>0\:.
\ees
Let now $\delta:=1-\frac{2}{p}\in(0,1)$, and let $0<\eta<p$ be arbitrary. We choose some $\epsilon$ satisfying
\bes
0\;<\;\epsilon\;<\;\frac{2\,\eta}{p(p-\eta)}\;\equiv\;\delta-1+\frac{2}{p-\eta}\:.
\ees
We have observed that $\,\pro\bA=\text{O}(|z|^{\delta-\epsilon})$, hence $\,
\pro\bA=\text{o}\big(|z|^{1-\frac{2}{p-\eta}}\big)$\,,
and in particular, we find
\be\label{pina}
|z|^{-1}\pi_{\bn}\bA\;\in\;L^{p-\eta}(\di)\:\:\:\:\qquad\forall\:\:\eta>0\:.
\ee
This fact shall come helpful in the sequel.

\medskip
When $\theta_0=1$, one directly deduces from the standard Calderon-Zygmund theorem applied to (\ref{delphi}) that $\nabla^2\bp\in L^p$. In that case, $\text{e}^{\la}$ is bounded from above and below. Owing to the identity
\bes
\nabla\bn\;=\;\star\,\Big[\big(\pro\nabla\bAe_1\big)\wedge\bAe_{2}\,+\,\bAe_{1}\wedge\big(\pro\nabla\bAe_2\big)\Big]\:,
\ees
we have
\be\label{modulegradn}
\big|\nabla\bn\big|\;=\;\text{e}^{-\la}\big|\pro\nabla^2\bp\big|\:,
\ee
thereby yielding that $\nabla\bn\in L^p$. When now $\theta_0\ge2$, we must proceed slightly differently to obtain analogous results. From (\ref{lambdas}), we know that $|z|\nabla\la$ is bounded across the unit-disk. We may thus apply Proposition \ref{CZpondere}-(ii) to (\ref{delphi}) with the weight $|\mu|=\text{e}^{\la}$ and $a=\theta_0-1$. The required hypothesis (\ref{hypw2}) is fullfilled, and we so obtain
\be\label{locexdelphi}
\nabla^2\bp\;=\;\theta_0\,(1-\theta_0)\left(\begin{array}{cc}-\,\Re&\Im\\[.5ex]\Im&\Re\end{array}\right)\big(\vec{A}\;{z}^{\theta_0-2}\big)\,+\,\text{e}^{\la}\vec{Q}\:,
\ee
where $\vec{A}$ is as in (\ref{locexphi}), and $\vec{Q}$ lies in $\mathbb{R}^4\otimes L^{p-\epsilon}(\di,\R^m)$ for every $\epsilon>0$. The exponent $p>2$ is the same as in (\ref{gradSinLp}).\\
Since $\,\text{e}^{\la}\simeq|z|^{\theta_0-1}$, we obtain from (\ref{locexdelphi}) that
\bes
\text{e}^{-\la}\big|\pro\nabla^2\bP\big|\;\lesssim\;|z|^{-1}|\pro\vec{A}|\,+\,|\pro\,\vec{Q}|\:.
\ees
According to (\ref{pina}), the first summand on the right-hand side of the latter belongs to $L^{p-\eta}\,$ for all $\eta>0$. Moreover, we have seen that $\pro\bQ$ lies in $L^{p-\epsilon}$ for all $\epsilon>0$. Whence, it follows that $\,\text{e}^{-\la}\pro\nabla^2\bP\,$ is itself an element of $L^{p-\epsilon}$ for all $\epsilon>0$. Brought into (\ref{modulegradn}), this information implies that
\be\label{regn}
\nabla\bn\;\in\;L^{p-\epsilon}(\di)\:,\qquad\quad\forall\:\:\epsilon>0\:.
\ee

In light of this new fact, we may now return to (\ref{sysSR}). In particular, recalling (\ref{regg}), we find 
\bes
\Delta S\;\equiv\;-\,\nabla(\star\,\bn)\cdot\big(\nabla^\perp\bR\,+\nabla\bG\big)\,-\,(\star\,\bn)\cdot\Delta\bG\:\:\in\:\:L^{q}(\di)\:,
\ees
with
\bes
\dfrac{1}{q}\:=\:\dfrac{1}{p}\,+\,\dfrac{1}{p-\epsilon}\:.
\ees
We attract the reader's attention on an important phenomenon occurring when $\theta_0=1$. In this case, if the aforementioned value of $q$ exceeds 2 (i.e. if $p>4$), then $\Delta S\notin L^q$, but rather only $\Delta S\in L^{2,\infty}$. This integrability ``barrier" stems from that of $\Delta\bG$, as given in (\ref{regg}). The same considerations apply of course with $\bR$ and $g$ in place of $S$ and $\bG$, respectively. \\[1ex]
Our findings so far may be summarized as follows:
\be\label{regsr}
\nabla S\;\,,\,\nabla\bR\:\,\in\:\,\left\{\begin{array}{lcl}W^{1,(2,\infty)}&,&\text{if}\:\:\:\theta_0=1\:\:\:\text{and}\:\:\:p>4\\[1ex]
W^{1,q}&,&\text{otherwise}.
\end{array}\right.
\ee
With the help of the Sobolev embedding theorem\footnote{we also use a result of Luc Tartar \cite{Ta} stating that $\,W^{1,(2,\infty)}\subset{BMO}$.}, we infer
\be\label{improvesr}
\nabla S\;\,,\,\nabla\bR\:\,\in\:\,\left\{\begin{array}{lcl}BMO&,&\text{if}\:\:\:\theta_0=1\:\:\:\text{and}\:\:\:p>4\\[1ex]
L^\infty&,&\text{if}\:\:\:\theta_0\ge2\:\:\:\text{and}\:\:\:p>4\\[1ex]
L^s&,&\text{if}\:\:\:\theta_0\ge1\:\:\:\text{and}\:\:\:p\le4\:,
\end{array}\right.
\ee
with
\bes
\dfrac{1}{s}\:=\:\dfrac{1}{q}\,-\,\dfrac{1}{2}\:\,=\:\,\dfrac{1}{p}\,+\,\dfrac{1}{p-\epsilon}\,-\,\dfrac{1}{2}\:\,<\,\:\dfrac{1}{p}\:.
\ees
Comparing (\ref{improvesr}) to (\ref{gradSinLp}), we see that the integrability has been improved. The process may thus be repeated until reaching that
\bes
\nabla S\:\,,\,\nabla\bR\:\,\in\:L^{b}(\di)\qquad\:\:\forall\:\:\:b<\infty
\ees
holds in all configurations. With the help of this newly found fact, we reapply Proposition \ref{CZpondere} so as to improve (\ref{regsr}) and (\ref{regn}) to
\bes
\nabla S\;\,,\,\nabla\bR\:\,\in\:\,\left\{\begin{array}{lcl}W^{1,(2,\infty)}(\di)&,&\text{if}\:\:\:\theta_0=1\\[1ex]
W^{1,b}(\di)&,&\text{if}\:\:\:\theta_0\ge2\:,\qquad\forall\:\:b<\infty
\end{array}\right.
\ees
and
\bes
\nabla\bn\,\in\,L^{b}(\di)\,\quad\quad\forall\:\:b<\infty\:.
\ees
The $\eps$-regularity in the form (\ref{barbibulle}) then yields a pointwise estimate for the Gauss map. Namely, in a neighborhood of the origin, 
\be\label{regn22}
|\nabla\bn|(x)\;\lesssim\;|x|^{-\epsilon}\qquad\forall\:\:\epsilon>0\:.
\ee

\medskip

This nearly completes the proof of Proposition \ref{Th1}. To fully end it, it remains to prove that for $\theta_0\ge2$, we may choose $\epsilon=0$ in (\ref{regn22}). This will be done in Section \ref{improv}, Remark \ref{remend}.

\subsection{Main results}
\subsubsection{Preparation}
We have recalled in the Introduction that it is shown in \cite{MS} that the conformal parameter satisfies
\bes
\la(z)\;=\;(\theta_0-1)\log|z|\,+\,u(z)
\ees
where the function $u$ belongs to $W^{2,1}$, and $\text{e}^{\pm u(0)}\ne0$. More precisely, from the Liouville equation, we know that
\be\label{deluu}
-\,\Delta u\;=\;\text{e}^{2\la}K\:,
\ee
with $K$ denoting the Gauss curvature. It is not hard to see that $\text{e}^{2\la}K$ inherits the regularity of $|\nabla\bn|^2$ (as it is made of products of terms of the type $\text{e}^{\la}\bh_{ij}$, each of which inherits the regularity of $|\nabla\bn|$). Owing to (\ref{deluu}), we see that $\,\nabla^2u\in\bigcap_{p<\infty}L^p$, and may thus write
\be\label{exla}
\text{e}^{-2\la(z)}\;=\;\big(\text{e}^{-2u(0)}+p(z)\big)|z|^{-2(\theta_0-1)}\:,
\ee
with $p(z)=\text{O}(|z|)$, and $p$ belongs to $C^{1,1-\epsilon}(\di)$ for all $\epsilon>0$.\\
Since
\bes
f(\bar{z})\;=\;{a_\mu}\bar{z}^{\,\mu}+f_0(\bar{z})\qquad\text{with}\quad\mu\ge-1\:,\quad a_\mu\in\C\setminus\{0\}\quad\text{and}\quad f_0\in C^{\infty}(\di)\:,
\ees
we use (\ref{exla}) along with the expansions (\ref{locexphi}) to obtain
\be\label{tresimp}
\text{e}^{-2\la}f\,\pz\bp\;=\;\pzb\bF_\mu+\vec{J}\:,
\ee
where
\be\label{defmu}
\bF_{\mu}(\bar{z})\;:=\;\dfrac{1}{2}\,\bA\;a_\mu\,\theta_0\,\text{e}^{-2u(0)}\left\{\begin{array}{lcl}
2\log|z|&,&\mu=\theta_0-2\\[1ex]
\frac{1}{\mu+2-\theta_0}\,\overline{z}^{\,\mu+2-\theta_0}&,&\mu\ne\theta_0-2\:\,;
\end{array}\right.
\ee
while $\vec{J}$ is the $\C^m$-valued function
\begin{eqnarray}\label{defJ}
\vec{J}&:=&\text{e}^{-2\la}f\pz\bp-\pzb\bF_\mu\nonumber\\[1ex]
&\equiv&\overline{z}^{\,\mu+1-\theta_0}\Big[z^{1-\theta_0}\text{e}^{-2u}\pz\bp-\big(z^{1-\theta_0}\text{e}^{-2u}\pz\bp\big)(0)   \Big]\:,
\end{eqnarray}
satisfying
\be\label{asympJ}
\vec{J}\;=\;\text{O}(|z|^{\mu+2-\theta_0})\:.
\ee
Converting (\ref{a1})-(\ref{a2}) into its complex form is easily done, namely
\begin{eqnarray*}
i\,\pzb\bL&=&\pzb\bH\,-\,3\,\pro\pzb\bH\,+\,i\star\big(\pzb\bn\wedge\bH\big)\,-\,\text{e}^{-2\la}f\,\pz\bp\,-\,2\,\vec{\beta}_0\,\pzb\log|z|\nonumber\\[1ex]
&=&-\,2\,\pzb\bH+3\,\pi_T\pzb\bH-i\star\big(\bn\wedge\pi_T\pzb\bH\big)-\text{e}^{-2\la}f\,\pz\bp-2\,\vec{\beta}_0\,\pzb\log|z|\:,
\end{eqnarray*}
where we have used the fact that $\bH$ is a normal vector\footnote{recall also that $\pi_T$ denotes projection onto the tangent bundle, namely $\pi_T:=\text{id}-\pro\,$.}.
Bringing (\ref{tresimp}) into the latter and rearranging the terms yields the identity
\bes
\pzb\big(i\bL+2\bH+2\vec{\beta}_0\log|z|+\bF_\mu\big)\;=\;3\,\pi_T\pzb\bH-i\star\big(\bn\wedge\pi_T\pzb\bH\big)-\bJ\:.
\ees
Using the fact that $\bH$ is normal, a simple computation reveals that
\begin{eqnarray}\label{unc1}
\pi_T\pzb\bH&=&-\sum_{j=1,2}\big(\bH\cdot\pro\pzb\vec{e}_{j}\big)\,\vec{e}_{j}\nonumber\\
&=&-\,(\bH\cdot\bH)\,\pzb\bp\,-\,(\bH\cdot{\bH^*_0})\,\pz\bp\:.
\end{eqnarray}
From this and the elementary identities
\bes
\star\,(\bn\wedge\pz\bp)\;=\;i\,\pz\bp\qquad\text{and}\qquad \star(\bn\wedge\pzb\bp)\;=\;-\,i\,\pzb\bp\:,
\ees
we obtain
\be\label{unc2}
\star\,(\bn\wedge\pi_T\pzb\bH)\;=\;i\,(\bH\cdot\bH)\,\pzb\bp\,-\,i\,(\bH\cdot{\bH^*_0})\,\pz\bp\:.
\ee
Altogether, (\ref{unc1}) and (\ref{unc2}) brought into (\ref{cayest2}) give
\begin{eqnarray}\label{cayest2}
\pzb\big(i\bL+2\bH+2\vec{\beta}_0\log|z|+\bF_\mu\big)&=&-\,2\,(\bH\cdot\bH)\,\pzb\bp\,-\,4\,(\bH\cdot{\bH^*_0})\,\pz\bp\,-\,\bJ\nonumber\\[1ex]
&=:&2\,\vec{q}\:.
\end{eqnarray}
Note that
\be\label{estq}
|\vec{q}|\;\lesssim\;|\nabla\bn||\bH|+|\bJ|\;\lesssim\;\text{O}\big(|z|^{-\epsilon}|\bH|+|z|^{\mu+2-\theta_0}\big)\:,
\ee
where we have used (\ref{regn22}) and (\ref{asympJ}).

\begin{Lma}\label{lemcrux}
Suppose that for some integer $k\in\{1,\ldots,\theta_0\}$, there holds locally around the origin
\be\label{Hhyp}
\bH\;=\;\text{O}(|z|^{k-\theta_0-\epsilon})\qquad\forall\:\epsilon>0\:,
\ee
and set
\bes
b\;:=\;\min\{k\,,\mu+2\}\:.
\ees
then we have
\be\label{correp}
\dfrac{i}{2}\,\bL+\bH+\vec{\beta}_0\log|z|+\dfrac{1}{2}\bF_\mu\;=\;\bE\,-\vec{T}\:.
\ee
The function $\bE$ is meromorphic with a pole at the origin of order
\bes
a\in\big\{\max\{0\,,\theta_0-\mu-2\},\ldots,\theta_0-b\big\}\:.
\ees
 Moreover
\bes
\pzb\bT\;=\;\vec{q}\quad\text{on}\:\:\di\setminus\{0\}\quad,\quad \bT\;=\;\text{O}(|z|^{1+b-\theta_0-\epsilon})\quad\forall\:\epsilon>0\:.
\ees
The function $\bT$ is unique up to addition of meromorphic summands.
\end{Lma}
\noindent
$\textbf{Proof.\:}$ 
Suppose that for some integer $k\in\{1,\ldots,\theta_0\}$ there holds
\bes
\bH\;=\;\text{O}(|z|^{k-\theta_0-\epsilon})\qquad\forall\:\epsilon>0\:.
\ees
We have seen that $|\nabla\bn|\lesssim|z|^{-\epsilon}$ for all $\epsilon>0$. Hence, since $|\bH||\pzb\bp|$ and $|\bH_0||\pz\bp|$ are controlled by $|\nabla\bn|$, it follows from  (\ref{asympJ}) and (\ref{cayest2}) that
\be\label{lem2q}
\vec{q}\;=\;\text{O}\big(|z|^{k-\theta_0-\epsilon}+|z|^{\mu+2-\theta_0}\big)\,=\;\text{O}(|z|^{b-\theta_0-\epsilon})\:.
\ee
We consider any $\vec{w}$ satisfying
\bes
\pzb\vec{w}\;=\;2\,{z}^{\theta_0-b}\vec{q}\quad\text{on}\:\:\di\:.
\ees
Per (\ref{lem2q}), $\vec{w}$ is $C^{0,1-\epsilon}$-H\"older continuous for any $\epsilon>0$. For notional convenience, let
\bes
\vec{W}\;:=\;\dfrac{i}{2}\,\bL+\bH+\vec{\beta}_0\log|z|+\dfrac{1}{2}\,\bF_\mu\:.
\ees
From (\ref{cayest2}), there holds
\be\label{enfinoolem}
\pzb\Big[{z}^{\theta_0-b}\,\vec{W}-\vec{w}\Big]\;=\;0\qquad\text{on}\:\:\di\setminus\{0\}\:.
\ee
We will extend this equation to all of the unit disk $\di$. To do so, it suffices to show that the function to which the operator $\pzb$ is applied on the left-hand side of the equation (\ref{enfinoolem}) lies in $L^2$. Since $\vec{w}$ is H\"older continuous, while $\bH$ satisfies (\ref{Hhyp}), the definition of $b$ guarantees that there only remains to verify that $|z|^{\theta_0-b}\bL$ lies in the space $L^2$. Exactly as we derived (\ref{xx2}) from (\ref{es2}), we infer here that 
$|z|^{\theta_0+1-b}\nabla\bH\in\bigcap_{p<\infty}L^p$, and then per (\ref{a1})-(\ref{a2}) that 
\begin{eqnarray}\label{gradLlem}
|z|^{\theta_0+1-b}|\nabla\bL|&\lesssim&|z|^{\theta_0+1-b}|\nabla\bH|+|z|^{\theta_0+1-b}|\nabla\bn||\bH|+|z|^{\mu+2-b}+|z|^{\theta_0-b}\nonumber\\[1ex]
&\in&\bigcap_{p<\infty}L^p(\di)\:,
\end{eqnarray}
from which we obtain, calling upon the Hardy-Sobolev inequality, that $|z|^{\theta_0-b}\bL\in L^2$. Accordingly, equation (\ref{enfinoolem}) holds on the unit disk. Whence,
\bes
\vec{W}\;=\;\vec{P}\,-\,{z}^{b-\theta_0}\vec{w}\:,
\ees
where $\vec{P}$ is meromorphic with a pole at the origin of order at most $(\theta_0-b)$. 
Putting in the latter
\bes
\bE\;:=\;\vec{P}+z^{b-\theta_0}\vec{w}(0)\qquad\text{and}\qquad \bT\;:=\;\big(\vec{w}-\vec{w}(0)\big)z^{b-\theta_0}
\ees
gives the desired representation (\ref{correp}). Moreover, there holds
\bes
\pzb\bT\;=\;\vec{q}\quad\text{on}\:\:\di\setminus\{0\}\quad,\quad \bT\;=\;\text{O}(|z|^{1+b-\theta_0-\epsilon})\quad\forall\:\epsilon>0\:.
\ees
The function $\vec{w}$ is clearly unique up to addition of meromorphic terms. The same is also true for $\bT$. Should the ``first" found $\bT$ happen to contain a meromorphic summand, it will necessarily be of order $|z|^{1+b-\theta_0}$ and could thus safely be fed into $\bE$ without affecting the desired statement.\\

The order $a$ of the pole of $\bE$ at the origin cannot possibly be strictly less than $\max\{0,\theta_0-\mu-2\}$. Indeed, it can certainly not be negative by definition. On the other hand, $a$ cannot either be strictly smaller than $\theta_0-\mu-2$. Indeed, if $\,0\le a\le \theta_0-\mu-3$, the representation (\ref{correp}) shows that
\bes
\bH+\vec{\beta}_0\log|z|\;=\;-\,\dfrac{1}{2}\,\Re\big(\bF_\mu\big)+\text{O}(|z|^{3+\mu-\theta_0-\epsilon})
\ees
In this configuration, $\theta_0\ge\mu+3$. Hence, $|\bF_\mu|\simeq|z|^{\mu+2-\theta_0}$. Accordingly, the latter shows that $\bH$ is dominated by $\frac{1}{2}\Re(\bF_\mu)$. This is however impossible. Indeed, $\bH$ is a normal vector, whereas $\bF_\mu$, being a multiple of the vector $\bA$ appearing in (\ref{delphi}), is a tangential vector near the origin. Hence the contradiction.\\[-3.5ex]

$\hfill\blacksquare$ \\

We now come to a central result in our study.
\begin{Prop}\label{Lemmeinter}
There exists a unique function $\bT$ containing no monomial of $z$, satisfying
\be\label{condernst}
\pzb\vec{T}\;=\;\vec{q}\quad\text{on}\:\:\di\setminus\{0\}
\qquad\text{and}\qquad \vec{T}\;=\;\text{O}\big(|z|^{2-\theta_0-\epsilon}\big)\quad\forall\:\epsilon>0\:,
\ee
and such that locally around the singularity, there holds
\be\label{thexphinter}
\dfrac{i}{2}\,\bL+\bH+\vec{\beta}_0\log|z|+\dfrac{1}{2}\,\bF_\mu\;=\;\vec{E}-\vec{T}\:,
\ee
where $\vec{\beta}_0$ is the residue defined in (\ref{residudu}), while the function $\vec{E}$ is holomorphic with possibly a pole at the origin of order at most $(\theta_0-1)$.\\[1ex]
Let $\al\ge0$ denote the order of the pole of the meromorphic function $\bE$ at the origin. Then
\begin{itemize}
\item[(i)] $\:\:\al\in\big\{\max\{0\,,\theta_0-\mu-2\},\ldots,\theta_0-1\big\}$\,;
\item[(ii)] the functions $\bE$ and $\bT$ may be adjusted to satisfy
\bes
\bE-\bT\;=\;\bE_\al z^{-\al}-\bQ\:,
\ees
for some nonzero constant $\bE_\al\in\C^m$, and with
\be\label{ilfo}
\pzb\bQ\;=\;\vec{q}\quad\text{on}\:\:\di\setminus\{0\}\quad,\quad \bQ\;=\;\text{O}(|z|^{1-a-\epsilon})\quad\forall\:\epsilon>0\:.
\ee
\end{itemize}
\end{Prop}
$\textbf{Proof.\:}$ 
We have seen in Proposition \ref{Th1} that $|\bH|\lesssim\text{e}^{-\la}|\nabla\bn|=\text{O}(|z|^{1-\theta_0-\epsilon})$ for all $\epsilon>0$. Since $\mu\ge-1$ always, we may thus choose $k=1=b$ in Lemma \ref{lemcrux}, thereby yielding both (\ref{thexphinter}) and (\ref{condernst}).\\

We next establish item (i) of the second part of the announced statement. We have already seen in the first part of the proposition that $a\le\theta_0-1$. Remains thus to establish that $a\ge\max\{0,\theta_0-\mu-2\}$.  This is trivially true for $\theta_0=1$, since $a\ge0$ by definition and $\mu\ge-1$ by hypothesis. We shall thus henceforth consider only the case $\theta_0\ge2$. For the sake of brevity, the claim is proved only for $\mu\in\{-1,0\}$. All other cases are obtained {\it mutatis mutandis}.

\paragraph{Case $\mu=-1$.} We need to show that $a\ge\max\{0,\theta_0-1\}=\theta_0-1$. Assume on the contrary that $0\le a\le \theta_0-2$. By hypothesis, 
\bes
\dfrac{i}{2}\,\bL+\bH+\vec{\beta}_0\log|z|+\dfrac{1}{2}\,\bF_\mu\;=\;\vec{E}-\vec{T}\:,
\ees
where $\bE$ is meromorphic with a pole of order $a$, and $\bT=\text{O}(|z|^{2-\theta_0-\epsilon})$. In particular, we find
\bes
\bH+\vec{\beta}_0\log|z|+\dfrac{1}{2}\,\Re\big(\bF_\mu)\;=\;\text{O}(|z|^{2-\theta_0-\epsilon})\:,
\ees
so that
\bes
\bH+\vec{\beta}_0\log|z|\;=\;-\,\dfrac{1}{2}\,\Re\big(\bF_\mu)+\text{O}(|z|^{2-\theta_0-\epsilon})\:.
\ees
ince $|\vec{F}_\mu|\simeq|z|^{\mu+2-\theta_0}$, it follows that $\bH$ is dominated $\Re(\vec{F}_\mu)$. This is however not possible, for $\bH$ is a normal vector, while $\bF_\mu$ is a tangent vector in a neighborhood of the origin.  

\paragraph{Case $\mu=0$.} We need to show that $a\ge\theta_0-2$. Assume on the contrary that $0\le a\le \theta_0-3$. By hypothesis, 
\be\label{rep1}
\dfrac{i}{2}\,\bL+\bH+\vec{\beta}_0\log|z|+\dfrac{1}{2}\,\bF_\mu\;=\;\vec{E}-\vec{T}\:,
\ee
where $\bE$ is meromorphic with a pole of order $a$, and $\bT=\text{O}(|z|^{2-\theta_0-\epsilon})$. In particular, we find
\bes
\bH+\vec{\beta}_0\log|z|+\dfrac{1}{2}\,\Re\big(\bF_\mu)\;=\;\text{O}(|z|^{2-\theta_0-\epsilon})\:,
\ees
so that
\bes
\bH\;=\;\text{O}(|z|^{2-\theta_0-\epsilon})\:.
\ees
Calling upon Lemma \ref{lemcrux} with $k=2$ gives us the representation
\be\label{rep2}
\dfrac{i}{2}\,\bL+\bH+\vec{\beta}_0\log|z|+\dfrac{1}{2}\,\bF_\mu\;=\;\vec{E}^1-\vec{T}^1\:,
\ee
where $\bE^1$ is meromorphic with a pole at the origin of order at most $(\theta_0-2)$, while $\bT^1=\text{O}(|z|^{3-\theta_0-\epsilon})$. Equating the identities (\ref{rep1}) and (\ref{rep2}) yields
\bes
\bT\;=\;\bE-\bE^1+\bT^1\:.
\ees
If $\bE^1$ has a pole of order exactly $(\theta_0-2)$, then it is not difficult to see that it dominates the right-hand side of the latter. Writing $\bE^1=\bE^1_0z^{2-\theta_0}+\text{O}(|z|^{3-\theta_0})$ for a nonzero constant $\bE^1_0\in\C^m$, we would whence find
\bes
\bT\;=\;\bE^1_0z^{2-\theta_0}+\,\text{O}(|z|^{3-\theta_0-\epsilon})\:.
\ees
This is impossible, since $\bT$ was chosen to contain no monomial of $z$. On the other hand, if $\bE^1$ has a pole of order strictly less than $(\theta_0-2)$, the representation (\ref{rep2}) gives
\bes
\bH+\vec{\beta}_0\log|z|+\dfrac{1}{2}\,\Re\big(\bF_\mu)\;=\;\text{O}(|z|^{3-\theta_0-\epsilon})\:.
\ees
As $\theta_0\ge3$ by hypothesis, $\Re(\bF_\mu)\simeq|z|^{2-\theta_0}$ dominates $\bH$, which, as explained in the previous case, is impossible. \\

Next, we establish claim (ii) of the lemma's statement. We let $\,\theta_0-1\ge a\ge\max\{0,\theta_0-\mu-2\}$. In particular, we are only concerned with cases for which the condition $\mu\ge\theta_0-a-2\,$ holds. Again for the sake of brevity, we only study in details the cases $a=\theta_0-1$ and $a=\theta_0-2$. All other cases are obtained {\it mutatis mutandis}.

\paragraph{Case $\al=\theta_0-1\,$.} 
We can write locally
\bes
\bE\;=\;\bE_{\theta_0-1}z^{1-\theta_0}+\,\bE^0\:,
\ees
where $\bE_{\theta_0-1}\in\C^m$ is constant, and $\bE^0$ is a meromorphic function with a pole at the origin of order at most $(\theta_0-2)$, i.e. $|\bE^0|\lesssim|z|^{2-\theta_0}$. We may then define $\bQ:=\bT-\bE^0$ satisfying $\pzb\vec{Q}=\pzb\bT$ on $\di\setminus\{0\}$, while $\bQ$ and $\bT$ have the same asymptotic behavior at the origin.
\paragraph{Case $\al=\theta_0-2\,$.}
In this case, without loss of generality, $\theta_0\ge2$ and $\mu\ge0$. 
As $\bE$ is meromorphic with a pole of order $(\theta_0-2)$ at the origin, we have $|\bE|\simeq|z|^{2-\theta_0}$ near the origin. The second condition in (\ref{condernst}) put into (\ref{thexphinter}) shows that $\bH$ is controlled by $|z|^{2-\theta_0-\epsilon}$ for all $\epsilon>0$. Calling upon Lemma \ref{lemcrux} with $k=2$ (and using that $\mu\ge0$) gives us the representation
\be\label{rep3}
\dfrac{i}{2}\,\bL+\bH+\vec{\beta}_0\log|z|+\dfrac{1}{2}\,\bF_\mu\;=\;\vec{E}^1-\vec{T}^1\:,
\ee
where $\bE^1$ is meromorphic with a pole at the origin of order at most $(\theta_0-2)$, while $\bT^1=\text{O}(|z|^{3-\theta_0-\epsilon})$. Equating the identities (\ref{thexphinter}) and (\ref{rep3}) yields
\bes
\bT\;=\;\bE-\bE^1+\bT^1\;=\;\bE-\bE^1+\text{O}(|z|^{3-\theta_0-\epsilon})\:.
\ees
As $\bT$ does not contain any monomial of $z$, while $\bE$ is holomorphic, we see that $\bE$ and $\bE^1$ must agree to highest order, hence $\bE^1$ has a pole of order $(2-\theta_0)$:
\bes
\bE^1\;=\;\bE_{\theta_0-2}z^{2-\theta_0}+\bE^2\:,
\ees
for some nonzero constant $\bE_{\theta_0-2}\in\C^m$, and where $\bE^2$ is a holomorphic function whose growth at the origin is controlled by $|z|^{3-\theta_0}$. Accordingly, (\ref{rep3}) yields
\bes
\dfrac{i}{2}\,\bL+\bH+\vec{\beta}_0\log|z|+\dfrac{1}{2}\,\bF_\mu\;=\;\bE_{\theta_0-2}z^{2-\theta_0}-\bQ\:,
\ees
where $\bQ:=\bE^2-\bT^1$. Clearly, $\bQ$ satisfies the required (\ref{ilfo}).  \\[-3.5ex]

$\hfill\blacksquare$ \\

We view $\bE$ as a string of $m$ complex-valued functions $\{E_j\}_{j=1,\ldots,m}$, all of which are meromorphic and may have a pole at the origin of order at least $\max\{0,\theta_0-\mu-2\}$ and at most $(\theta_0-1)$. In particular, we define the $\mathbb{N}^m$-valued {\it second residue}
\be\label{residudubis} 
\vec{\gamma}\;=\;(\gamma_1, \ldots, \gamma_m)\qquad\text{with}\qquad\gamma_j\;:=\;-\,\dfrac{1}{2i\pi}\int_{\partial\di} d\log{E_j} \:.
\ee

\subsubsection{How the second residue $\vec{\gamma}$ controls the regularity}\label{improv}

We start by defining
\bes
\al\;:=\,\max_{1\leq j\leq m}\gamma_j\,\in\,\big\{\max\{0\,,\theta_0-\mu-2\},\ldots,\theta_0-1\big\}\:. 
\ees
Per Proposition \ref{Lemmeinter}, we may choose $\bE(z)=\bE_\al z^{-\al}$ for some constant vector $\bE_\al\in\C^m$. According to Proposition \ref{Lemmeinter}, there holds
\be\label{Qdec}
\bQ\;=\;\text{O}(|z|^{1-\al-\epsilon})\qquad\forall\:\epsilon>0\:.
\ee
Because $\bL$ is real-valued, (\ref{thexphinter}) yields
\bes
\bH+\vec{\beta}_0\log|z|+\dfrac{1}{2}\,\Re\big(\bF_\mu)\;=\;\Re\big(\bE-\vec{Q}\big)\:.
\ees
The function $\bF_\mu$ displays different behaviors depending upon whether $\theta_0=\mu+2$ or not (cf. (\ref{defmu})). It is convenient to highlight this fact by recasting the latter in the form
\be\label{luti}
\bH+\vec{\gamma}_0\log|z|\;=\;\Re\big(\bF-\vec{Q}\big)\:,
\ee
where
\be\label{modifres}
\vec{\gamma}_0\;:=\;\vec{\beta}_0+\dfrac{1}{2}\,\delta_{\theta_0,\mu+2}\,\text{e}^{-2u(0)}\theta_0\,\Re(a_\mu\bA)
\ee
and
\be\label{Fdec}
\bF\;:=\;\bE\,-\,\dfrac{1}{2}(1-\delta_{\theta_0,\mu+2})\bF_\mu\;=\;\text{O}(|z|^{-a})\:.
\ee
Note that $\bF$ is a power function of $z$ and $\bar{z}$. \\

We define next the two-component vector-field $\,\vec{U}:=\big(\Re\,,\Im\big)(\text{e}^{\la}\vec{Q}\,)$. As $\pzb\bQ=\vec{q}$, there holds\footnote{although the equation for $\bQ$ holds only on $\di\setminus\{0\}$, the system for $\vec{U}$ may easily be extended to the whole unit disk $\di$ owing to the fact that $\vec{U}=\text{O}(|z|^{1-\epsilon})\in L^\infty$.}
\be\label{divcurlu}
\left\{\begin{array}{rcl}
\text{div}\;\vec{U}&=&\nabla\la\cdot\vec{U}\,-\,2\,\text{e}^{\la}\,\Re({\vec{q}})\\[1ex]
\text{curl}\;\vec{U}&=&\nabla^\perp\la\cdot\vec{U}\,+\,2\,\text{e}^{\la}\,\Im({\vec{q}})\:.
\end{array}\right.
\ee
Because $|\nabla\la|\lesssim|z|^{-1}$, the estimates (\ref{estq}) and (\ref{Qdec}) give
\begin{eqnarray}\label{estimatU}
\big|\text{div}\;\vec{U}\big|+\big|\text{curl}\;\vec{U}\big|&\lesssim&|z|^{-1}|\vec{U}|+\text{e}^{\la}|z|^{-\epsilon}|\bH|+\text{e}^{\la}|z|^{\mu+2-\theta_0}\nonumber\\[1ex]
&\lesssim&|z|^{\theta_0-2-\epsilon}|\bQ|+|z|^{\theta_0-1-\epsilon}|\bF|+|z|^{\mu+1}\nonumber\\[1ex]
&\lesssim&|z|^{\theta_0-1-\al-\epsilon}\qquad\forall\:\epsilon>0\:,
\end{eqnarray}
where we have used that $\,a\ge\theta_0-\mu-2$.\\
With the help of a simple Hodge decomposition, (\ref{estimatU}) along with the fact that $|\vec{U}|\simeq\text{e}^{\la}|\bQ|=\text{O}(|z|^{\theta_0-a-\epsilon})$ yields
\bes
\big|\nabla\big(\text{e}^{\la}\bQ\big)\big|\;\simeq\;|\nabla\vec{U}|\;\lesssim\;|z|^{\theta_0-1-a-\epsilon}\qquad\forall\:\epsilon>0\:.
\ees
Again since $|\nabla\la|\lesssim|z|^{-1}$, the latter shows that 
\be\label{yaka1}
|\nabla\bQ|\;\lesssim\;|z|^{-a-\epsilon}\qquad\forall\:\epsilon>0\:.
\ee
Since $\bF$ is a power function of order $(-a)$, there holds
\be\label{estimatE}
|\nabla\bF|\;\lesssim\;|z|^{-1-a}\:.
\ee
Putting (\ref{yaka1}) and (\ref{estimatE}) into (\ref{luti}) then yields
\be\label{delH1}
|\nabla\bH|\:\lesssim\:|z|^{-1-\al}\:.
\ee
As $\al\le\theta_0-1$, we thus find $\text{e}^{\la}\nabla\bH\in L^{2,\infty}$. It is proved in \cite{BR3} (Section A.2.1) that the $\Lambda^{m-2}(\mathbb{S}^{m-1})$-valued Gauss map $\bn$ satisfies a perturbed harmonic map equation:
\be\label{prof444}
\Delta\bn\,-\,2\,\text{e}^{2\la}K\,\bn\;=\;2\star\!\big(\nabla^\perp\bp\wedge\nabla\bH\big)\,-\,2\star\text{e}^{2\la}\vec{h}_{12}\wedge(\vec{h}_{11}-\vec{h}_{22})\:,
\ee
where $K$ is the Gauss curvature. Whence,
\be\label{ghuy}
|\Delta\bn|\:\lesssim\:\text{e}^{\la}|\nabla\bH|\,+\,|\nabla\bn|^2\:\lesssim\:|z|^{\theta_0-2-a}\:\in\:L^{2,\infty}\:,
\ee
since $a\le\theta_0-1$. Accordingly, $\nabla^2\bn\in L^{2,\infty}$, and in particular $\nabla\bn\in BMO$.\\[1ex]

We have seen in the Introduction that the conformal parameter satisfies
\be\label{mulsve10}
\la\;=\;(\theta_0-1)\log|z|\,+u\:,
\ee
where the function $u$ belongs to $W^{2,1}$. More precisely, from the Liouville equation, we know that
\bes
-\,\Delta u\;=\;\text{e}^{2\la}K\:,
\ees
with $K$ denoting the Gauss curvature. As previously explained, $\text{e}^{2\la}K$ inherits the regularity of $|\nabla\bn|^2$ (as it is made of products of terms of the type $\text{e}^{\la}\bh_{ij}$, each of which inherits the regularity of $|\nabla\bn|$). Owing to (\ref{ghuy}), we thus have that $\,\nabla^2u\in\bigcap_{p<\infty}L^p$, and in particular that $\nabla u$ is H\"older continuous. Hence, (\ref{mulsve10}) shows that
\be\label{diffla}
|\nabla\la|\:\lesssim\:|z|^{-1}\:.
\ee
Furthermore, we may write
\be\label{decla1}
2\,\text{e}^{2\la}\;=\;\big(T_1+R_1\big)|z|^{2(\theta_0-1)}\:,
\ee
where $T_1$ is the first-order Taylor polynomial expansion of $2\text{e}^{2u}\in C^{1,1-\epsilon}$ (for all $\epsilon>0$) near the origin, and $R_1$ is the corresponding remainder. Hence
\be\label{R1info}
\nabla^jR_1\,=\,\text{O}(|z|^{2-j-\epsilon})\:,\quad j\in\{0,1\}\:,\:\:\:\forall\:\epsilon>0\:.
\ee
With the help of (\ref{luti}), we write
\bes
\Delta\bp\;\equiv\;2\,\text{e}^{2\la}\bH\;=\;\Delta\bp_0\,+\,\Delta\bp_1\:,
\ees
where
\bes
\left\{\begin{array}{lcl}\Delta\bp_0&=&T_1|z|^{2(\theta_0-1)}\Re\big(\bF-\vec{\gamma}_0\log|z|\big)\\[1.5ex]
\Delta\bp_1&=&-\,2\,\text{e}^{2\la}\Re(\bQ)+|z|^{2(\theta_0-1)}R_1\Re\big(\bF-\vec{\gamma}_0\log|z|\big)\:.\end{array}\right.
\ees
Since $T_1$ and $\bF$ are power functions, we easily obtain via solving explicitly and handling the remainder with Proposition \ref{CZcoro} that
\bes
\bp_0\;=\;\Re(\vec{P}_0)+\,C_1|z|^{2\theta_0}\Re(\bF)+\bC\,|z|^{2\theta_0}\big(\log|z|^{\theta_0}-1\big)+\vec{\xi}_0\:,
\ees
where $\vec{P}_0$ is a $\C^m$-valued holomorphic polynomial of degree at most $(2\theta_0-a)$, and 
\bes
C_1\;:=\;\dfrac{\text{e}^{2u(0)}}{2\theta_0}\qquad\text{and}\qquad \bC\;:=\;\dfrac{\text{e}^{2u(0)}}{2\theta_0^2}\,\vec{\gamma}_0\:.
\ees
The remainder $\vec{\xi}_0$ satisfies
\bes
\nabla^j\vec{\xi}_0\;=\;\text{O}(|z|^{2\theta_0-a+1-j-\epsilon})\qquad\forall\,\:j\in\{0,\ldots,2\}\:,\:\:\forall\:\,\epsilon>0\:,
\ees
and
\bes
|z|^{2+a-2\theta_0}\nabla^3\vec{\xi}_0\,\in\bigcap_{p<\infty}L^p\:.
\ees
To obtain information on $\bp_1$, we differentiate once its partial differential equation in each coordinate $x_1$ and $x_2$, and apply Proposition \ref{CZcoro} to the yield, using (\ref{R1info}), the fact that $\text{e}^{\la}\bQ\in\bigcap_{p<\infty}W^{1,p}$, that $\nabla\la=\text{O}(|z|^{-1})$, and the fact that $\bF$ is a power function. Without much effort, it ensues that we can write
\bes
\bp_1\;=\;\Re(\vec{P}_1)+\vec{\xi}_1\:,
\ees
where $\vec{P}_1$ is a $\C^m$-valued holomorphic polynomial of degree at most $(2\theta_0-a)$, and the $\R^m$-valued function $\vec{\xi}_1$ satisfies
\bes
\nabla^j\vec{\xi}_1\;=\;\text{O}(|z|^{2\theta_0-a+1-j-\epsilon})\qquad\forall\,\:j\in\{0,\ldots,2\}\:,\:\:\forall\:\,\epsilon>0\:,
\ees
and
\bes
|z|^{2+a-2\theta_0}\nabla^3\vec{\xi}_1\,\in\bigcap_{p<\infty}L^p\:.
\ees
Comparing $\bp_0+\bp_1$  to the previously found expression (\ref{locexphi}), we deduce
\be\label{gary13}
\bp=\Re\big(\bA\,z^{\theta_0}+\bB_1z^{\theta_0+1}+C_{1}|z|^{2\theta_0}\bF\big)+\bC|z|^{2\theta_0}\big(\log|z|^{\theta_0}-1\big)+(\vec{\xi}_0+\vec{\xi}_1)\:,
\ee
where $\bB_1\in\C^m$ is constant, while $\bA$ is as in Proposition \ref{Th2}. \\
Note that
\be\label{2diffphi}
|\nabla^j\bp|\;=\;\text{O}(|z|^{\theta_0-j})\qquad\forall\:j\in\{0,\ldots,2\}\:.
\ee
\smallskip

Suppose next that $\,\theta_0-2\ge a\ge\theta_0-\mu-2$. Then (\ref{delH1}) gives
\be\label{delH2}
\text{e}^{\la}\nabla\bH\,\in\,L^\infty\:.
\ee
In turn brought into (\ref{ghuy}), the latter shows that
\be\label{diff2n}
\nabla^2\bn\,\in\bigcap_{p<\infty}L^p\:.
\ee
Accordingly, the function $u$ appearing in (\ref{mulsve10}) lies in $C^{2,1-\epsilon}(\di)$ for all $\epsilon>0$, whence
\bes
\nabla^2\la\;=\;\text{O}(|z|^{-2})\:.
\ees
When $a\le\theta_0-2$, we have that $|z|^{-1}\text{e}^{\la}|\bF|\simeq|z|^{\theta_0-2-\al}\in L^\infty$. Hence (\ref{Qdec}) and (\ref{luti}) yield
\be\label{hsurx}
|z|^{-1\,}\text{e}^{\la}\bH\;\in\;L^\infty\:.
\ee
We now need to improve the regularity of $\vec{q}$. Recall that
\bes
\vec{q}\;:=\;-\,|\bH|^2\pzb\bp\,-\,2\,(\bH\cdot\bH_0^*)\pz\bp\,-\,\dfrac{1}{2}\,\vec{J}\:,
\ees
where
\bes
\vec{J}\;:=\;\text{e}^{-2\la}f\pz\bp-\pzb\bF_\mu\;\equiv\;\overline{z}^{\,\mu+1-\theta_0}\Big[z^{1-\theta_0}\text{e}^{-2u}\pz\bp-\big(z^{1-\theta_0}\text{e}^{-2u}\pz\bp\big)(0)   \Big]\:,
\ees
where $u$ is the function appearing in (\ref{mulsve10}). \\
As we are studying the case $\theta_0-2\ge a\ge\theta_0-2-\mu$, we have automatically that $\mu\ge0$. Hence, using (\ref{2diffphi}) and the fact that $u\in C^1$, we see that
\bes
\nabla\bJ\;=\;\text{O}(|z|^{\mu+1-\theta_0})\:.
\ees
With (\ref{asympJ}) and (\ref{diffla}), the latter yields
\be\label{diffJ}
\nabla\big(\text{e}^{\la}\bJ\,\big)\;=\;\text{O}(|z|^{\mu})\,\in\,L^\infty\:.
\ee
On the other hand, since $\text{e}^{\la}\bH$ and $\text{e}^{\la}\bH_0$ inherit the regularity of $\nabla\bn$, we find
\begin{eqnarray}\label{lat2}
\big|\nabla\big(\text{e}^{\la}(\bH\cdot\bH_0^*)\pz\bp\big)\big|&\lesssim&|\text{e}^{\la}\bH||\nabla^2\bn|+|\nabla\bn|\big(\text{e}^{\la}|\nabla\bH|+|\nabla^2\bp||\bH|\big)\nonumber\\[1ex]
&\lesssim&|\nabla^2\bn|+\text{e}^{\la}|\nabla\bH|+|z|^{-1}\text{e}^{\la}|\bH|\nonumber\\[1ex]
&\in&\bigcap_{p<\infty}L^p\:,
\end{eqnarray}
where we have used successively, (\ref{diff2n}), (\ref{2diffphi}), (\ref{delH2}), and (\ref{hsurx}). Exactly in the same fashion, one verifies that
\bes
\text{e}^{\la}|\bH|^2\pzb\bp\;\in\bigcap_{p<\infty}W^{1,p}\:.
\ees
Together, the latter, (\ref{lat2}), and (\ref{diffJ}) brought into the definition of $\vec{q}$ show that
\be\label{diff1q}
\text{e}^{\la}\vec{q}\;\in\bigcap_{p<\infty}W^{1,p}\:.
\ee
We next return to the system (\ref{divcurlu}). Proceeding as in (\ref{estimatU}) with the information that $a\le\theta_0-2$, we infer that
\bes
\Big|\text{div}\big(|z|^{-1}\vec{U}\big)\Big|\,+\,\Big|\text{curl}\big(|z|^{-1}\vec{U}\big)\Big|\:\:\lesssim\:\:|z|^{\theta_0-2-\al-\epsilon}\;\lesssim\;|z|^{-\epsilon}\qquad\forall\:\epsilon>0\:,
\ees
so that $\,|z|^{-1}\text{e}^{\la}\bQ\equiv|z|^{-1}\vec{U}$ is an element of $W^{1,p}$ for all finite $p$. By a similar token, using (\ref{diff1q}), it is not difficult to see that
\bes
|\nabla^2\vec{U}|\;\lesssim\;\big|\nabla\big(|z|^{-1}\vec{U}\big)\big|+\big|\nabla(\text{e}^{\la}\vec{q})\big|\;\in\bigcap_{p<\infty}L^p\:.
\ees
Whence, $\,\text{e}^{\la}\vec{Q}\equiv\vec{U}\in\bigcap_{p<\infty}W^{2,p}$. \\
Using that $\nabla\la=\text{O}(|z|^{-1})$ now gives
\bes
|z|^{-1}\text{e}^{\la}|\nabla\bQ|\:\lesssim\:\big|\nabla\big(|z|^{-1}\text{e}^{\la}\bQ\big)\big|+|z|^{-2}\text{e}^{\la}|\bQ|\;\in\bigcap_{p<\infty}L^p\:,
\ees
where have used that $a\le\theta_0-2$ and $\bQ=\text{O}(|z|^{1-a-\epsilon})$ for all $\epsilon<0$. In particular, owing to (\ref{luti}) and (\ref{Fdec}), there holds
\bes
|z|^{-1}\text{e}^{\la}\big|\nabla(\bH+\vec{\gamma}_0\log|z|-\Re(\bF))\big|\:\lesssim\:|z|^{-1}\text{e}^{\la}\big|\nabla\bQ|\;\in\bigcap_{p<\infty}L^p\:.
\ees
Analogously, using now additionally that $\nabla^2\la=\text{O}(|z|^{-2})$ yields
\bes
\text{e}^{\la}|\nabla^2\bQ|\:\lesssim\:|z|^{-2}|\vec{U}|+\big|\nabla\big(|z|^{-1}\vec{U}\big)\big|+|\nabla^2\vec{U}|\:.
\ees
As we have shown above, each of these terms lies in $L^p$ for all finite $p$. Accordingly, differentiating twice (\ref{luti}) yields 
 \be\label{diff2h}
\text{e}^{\la}\big|\nabla^2(\bH+\vec{\gamma}_0\log|z|-\Re(\bF))\big|\:\lesssim\:\text{e}^{\la}\big|\nabla^2\bQ|\;\in\bigcap_{p<\infty}L^p\:.
\ee
We have already pointed out that the function $u$ in (\ref{mulsve10}) lies in $C^{2,1-\epsilon}$ for all $\epsilon>0$, owing to the fact that $\bn\in W^{2,p}$ for all $p<\infty$. We may now replace (\ref{decla1}) by
\bes
2\,\text{e}^{2\la}\;=\;\big(T_2+R_2\big)|z|^{2(\theta_0-1)}\:,
\ees
where $T_2$ is the second-order Taylor polynomial expansion of $2\text{e}^{2u}$, and $R_2$ is the corresponding remainder. Hence
\be\label{R2info}
\nabla^jR_2\,=\,\text{O}(|z|^{3-j-\epsilon})\:,\quad j\in\{0,\ldots,2\}\:,\:\:\:\forall\:\epsilon>0\:.
\ee
As before, we decompose
\bes
\Delta\bp\;\equiv\;2\,\text{e}^{2\la}\bH\;=\;\Delta\bp_0\,+\,\Delta\bp_1\:,
\ees
with now
\bes
\left\{\begin{array}{lcl}\Delta\bp_0&=&T_2|z|^{2(\theta_0-1)}\Re\big(\bF-\vec{\gamma}_0\log|z|\big)\\[1.5ex]
\Delta\bp_1&=&-\,2\,\text{e}^{2\la}\Re(\bQ)+|z|^{2(\theta_0-1)}R_2\Re\big(\bF-\vec{\gamma}_0\log|z|\big)\:.\end{array}\right.
\ees
Since $T_2$ and $\bF$ are power functions, we easily obtain via solving explicitly and handling the remainder with Proposition \ref{CZcoro} that
\bes
\bp_0\;=\;\Re(\vec{P}_0)+\,C_2\,|z|^{2\theta_0}\Re(\bF)+\bC\,|z|^{2\theta_0}\big(\log|z|^{\theta_0}-1\big)+\vec{\xi}_0\:,
\ees
where $\vec{P}_0$ is a $\C^m$-valued holomorphic polynomial of degree at most $(2\theta_0-a)$, and 
\bes
C_2\;:=\;\dfrac{\text{e}^{2u(0)}}{4\theta_0}\qquad\text{and}\qquad \bC\;:=\;\dfrac{\text{e}^{2u(0)}}{2\theta_0^2}\,\vec{\gamma}_0\:.
\ees
The remainder $\vec{\xi}_0$ satisfies
\bes
\nabla^j\vec{\xi}_0\;=\;\text{O}(|z|^{2\theta_0-a+1-j-\epsilon})\qquad\forall\,\:j\in\{0,\ldots,3\}\:,\:\:\forall\:\,\epsilon>0\:,
\ees
and
\bes
|z|^{3+a-2\theta_0}\nabla^4\vec{\xi}_0\,\in\bigcap_{p<\infty}L^p\:.
\ees
To obtain information on $\bp_1$, we differentiate twice its partial differential equation in each coordinate $x_1$ and $x_2$, and apply Proposition \ref{CZcoro} to the yield, using (\ref{R2info}), the fact that $\text{e}^{\la}\bQ\in\bigcap_{p<\infty}W^{2,p}$, that $\nabla^2\la=\text{O}(|z|^{-2})$, and the fact that $\bF$ is a power function. Without much effort, it ensues that we can write
\bes
\bp_1\;=\;\Re(\vec{P}_1)+\vec{\xi}_1\:,
\ees
where $\vec{P}_1$ is a $\C^m$-valued holomorphic polynomial of degree at most $(2\theta_0-a)$, and the $\R^m$-valued function $\vec{\xi}_1$ satisfies
\bes
\nabla^j\vec{\xi}_1\;=\;\text{O}(|z|^{2\theta_0-a+1-j-\epsilon})\qquad\forall\,\:j\in\{0,\ldots,3\}\:,\:\:\forall\:\,\epsilon>0\:,
\ees
and
\bes
|z|^{3+a-2\theta_0}\nabla^4\vec{\xi}_1\,\in\bigcap_{p<\infty}L^p\:.
\ees
Comparing $\bp_0+\bp_1$  to the previously found expression (\ref{gary13}), we deduce
\begin{eqnarray}\label{gary14}
\bp&=&\Re\big(\bA\,z^{\theta_0}+\bB_1z^{\theta_0+1}+\bB_2z^{\theta_0+2}+C_{2}|z|^{2\theta_0}\bF\big)\nonumber\\[0ex]
&&\hspace{1.5cm}+\:\;\bC|z|^{2\theta_0}\big(\log|z|^{\theta_0}-1\big)+(\vec{\xi}_0+\vec{\xi}_1)\:,
\end{eqnarray}
where $\bA$ and $\bB_1$ are as in (\ref{gary13}), while $\bB_2\in\C^m$ is constant.\\
Note that
\be\label{3diffphi}
|\nabla^j\bp|\;=\;\text{O}(|z|^{\theta_0-j})\qquad\forall\:j\in\{0,\ldots,3\}\:.
\ee
\smallskip
Finally, we return to the equation (\ref{prof444}). Using the previously noted fact that $\text{e}^{\la}\bh_{ij}$ inherit the regularity of $\nabla\bn$, along with (\ref{diff2n}), (\ref{diff2h}), (\ref{3diffphi}), we now obtain
\begin{eqnarray*}
|\Delta\nabla\bn|&\lesssim&|\nabla^2\bn|\,+\,|\nabla\bn|^2|\nabla\bn|\,+\,|\nabla\bp|\,|\nabla^2\bH|\,+\,|\nabla^2\bp|\,|\nabla\bH|\nonumber\\[1ex]
&\lesssim&|\nabla^2\bn|\,+\,|\nabla\bn|^2|\nabla\bn|\,+\,\text{e}^{\la}\,|\nabla^2\bH|\,+\,|z|^{-1}\text{e}^{\la}\,|\nabla\bH|\nonumber\\[1ex]
&\simeq&|z|^{\theta_0-3-\al}+\,\text{terms in}\:\bigcap_{p<\infty}L^p\:.
\end{eqnarray*}
This shows that $\nabla^3\bn\in L^{2,\infty}$ if $\al=\theta_0-2$. On the other hand, if $a\le\theta_0-3$, we obtain that $\bn\in\bigcap_{p<\infty}W^{3,p}$. We may then start over again the above procedure gaining one order of decay at every step. The condition that $\,a\ge\max\{0,\theta_0-\mu-2\}$ is essential. It guarantees indeed that $\mu$ is large enough for a given integer $a$. In particular, it is possible to control the regularity of $\vec{J}$ appearing in the definition of $\vec{q}$. This is contrast with the case of Willmore immersions for which the multiplier function $f$ (and thus $\vec{J}$) are identically zero. In rough terms, this amounts to choosing $\mu=\infty$ in the present routine.\\

From the cases $a=\theta_0-1$ and $a=\theta_0-2$ treated above, a clear pattern emerges. Repeating finitely many times the steps performed above, one eventually reaches that
\be\label{item3A}
\nabla^{\theta_0-\al+1}\bn\;\in\,L^{2,\infty}\qquad\text{and thus}\qquad\nabla^{\theta_0-a}\bn\;\in\,BMO\:.
\ee
Furthermore, for all $\,j\in\{0,\ldots,\theta_0-\al\}$, there holds
\be\label{item2}
|z|^{a+j-1}\nabla^j\big(\bH+\vec{\gamma}_0\log|z|-\Re(\bF)\big)\;\in\bigcap_{p<\infty}L^p\:.
\ee
We also obtain a local expansion for the immersion, namely
\be\label{item3}
\bp=\Re\bigg(\bA\,z^{\theta_0}+\sum_{j=1}^{\theta_0-\al}\!\bB_j\,z^{\theta_0+j}+C_{\theta_0-\al}|z|^{2\theta_0}\bF\bigg)+\bC|z|^{2\theta_0}\big(\log|z|^{\theta_0}-1\big)+\vec{\xi}\:,
\ee
where $\bB_j\in\C^m$ are constant vectors, while $\bA$ is as in (\ref{thexphinter}). The constants $C_{\theta_0-a}$ and $\bC$ are
\bes
C_{\theta_0-a}\;:=\;\dfrac{\text{e}^{2u(0)}}{2\theta_0(\theta_0-a)}\qquad\text{and}\qquad \bC\;:=\;\dfrac{\text{e}^{2u(0)}}{2\theta_0^2}\,\vec{\gamma}_0\:.
\ees
The remainder $\vec{\xi}$ satisfies
\bes
\nabla^j\vec{\xi}\;=\;\text{O}(|z|^{2\theta_0-a+1-j-\epsilon})\qquad\forall\,\:j\in\{0,\ldots,\theta_0-a+1\}\:,\:\:\forall\:\,\epsilon>0\:,
\ees
and
\bes
|z|^{1-\theta_0}\nabla^{\theta_0-a+2}\vec{\xi}\,\in\bigcap_{p<\infty}L^p\:.
\ees
whence
\bes
\left\{\begin{array}{rclcl}
\nabla^{j}\bp&\!\!\!=\!\!\!&\text{O}\big(|z|^{\theta_0-j}\big)&,&j\in\{0,\ldots,\theta_0-a+1\}\\[1.25ex]
\nabla^{\theta_0+1}\bp&\!\!\!=\!\!\!&
\left\{\begin{array}{lcl}
\text{O}\big(|\log|z||\big)&\:,\:&\theta_0=1\\[1ex]
\text{O}(1)&\:,\:&\theta_0\ge2
\end{array}\right.
&,&\text{in the case}\:\:a=0\:.
\end{array}\right.
\ees

\smallskip
\noindent
Of course, when $\al>0$, the term $\vec{\xi}$ in (\ref{item3}) dominates the logarithmic term,  written here to indicate the presence and the influence of the (modified) first residue $\vec{\gamma}_0$ of which it is a multiple. Furthermore, the aforementioned information and the fact that $\bF$ is a power function of order $(-a)$ show that
\be\label{item37}
\bp\,\in\bigcap_{p<\infty}\left\{\begin{array}{lcl}
W^{\theta_0+2-a,p}&\:,\:&\theta_0\ge2\\[1ex]
W^{2,p}&\:,\:&\theta_0=1\\[1ex]
W^{3,p}&\:,\:&\theta_0=1\:\:,\:\:\vec{\gamma}_0=\vec{0}\:.
\end{array}\right.
\ee

\noindent
Using the definitions of $\bF$ and $\bE$ given at the beginning of this section, and using the fact that $a\ge\theta_0-\mu-2$, it is possible to reformulate (\ref{item2}) in the form
\be\label{item77}
|z|^{a+j-1}\nabla^j\big(\bH+\vec{\gamma_0}\log|z|-\Re(\bE_{a}z^{-a})\big)\in\bigcap_{p<\infty}L^p\qquad\forall\:j\in\{0,\ldots,\theta_0-a\}\,
\ee
where $\bE_{a}\in\C^m$ is the previously defined constant vector. 

\medskip

\begin{Rm}\label{remend}
At last we can complete the proof of Proposition \ref{Th1}, namely that $\nabla\bn$ is bounded across the unit disk when $\theta_0\ge2$. To see this, we first note that (\ref{item3}) yields
\bes
\pz\bp\;=\;\dfrac{\theta_0}{2}\,\bA\,z^{\theta_0-1}+\text{O}(|z|^{\theta_0})\:.
\ees
Since $\,\pro\nabla\bp\equiv0$, the latter gives $\,|\pi_{\bn}\bA|\;=\;\text{O}(|z|)$.
Differentiating twice (\ref{item3}) then reveals that
\bes
\big|\pi_{\bn}\nabla^2\bp\big|\:\lesssim\:|\pi_{\bn}\bA|\,|z|^{\theta_0-2}+\text{O}(|z|^{\theta_0-1})\;=\;\text{O}(|z|^{\theta_0-1})\:.
\ees
Combining this to (\ref{modulegradn}) gives thus that $\nabla\bn$ is bounded across the singularity. 
\end{Rm}

\subsubsection{When both residues vanish: improved regularity}\label{wbrv}

This last section is devoted to proving Theorem \ref{Th44}. We shall assume that the {\it modified first} residue $\vec{\gamma}_0$ and {\it second} residue $\vec{\gamma}$ defined respectively in (\ref{modifres}) and in (\ref{residudubis}) both vanish. Note that for the second residue to vanish, it is necessary to request that $\theta_0\le\mu+2$ (since the maximum component of $\vec{\gamma}$ is bounded from below by $\theta_0-\mu-2$). We will in time have to distinguish two cases depending upon whether $\theta_0<\mu+2$, in which case the immersion will be smooth, or whether $\theta_0=\mu+2$, in which case the immersion might be as little regular as $C^{2,\alpha}(\di)$, for all $\alpha<1$ (although $\nabla^{\theta_0+2}\bp$ always lies $L^p$ for all $p<\infty$).\\

Suppose that $\vec{\gamma}_0=\vec{0}=\vec{\gamma}$. Then the results of the previous section apply (with $a=0$). In particular, (\ref{item3A})-(\ref{item77}) give
\be\label{itemn1}
\nabla^{\theta_0}\bn\;\in BMO\:,
\ee
as well as
\bes
\left\{\begin{array}{rclcl}
\nabla^{j}\bp&\!\!\!=\!\!\!&\text{O}\big(|z|^{\theta_0-j}\big)&\:,\:&j\in\{0,\ldots,\theta_0\}\\[1ex]
|z|^{j-1}\nabla^j\bH&\!\!\!\in\!\!\!&\bigcap_{p<\infty}L^p&\:,\:&j\in\{1,\ldots,\theta_0\}\:.
\end{array}\right.
\ees
The latter yields
\bes
\big|\nabla^{\theta_0-1}(\nabla^\perp\bp\wedge\nabla\bH)\big|\:\:\lesssim\:\:\sum_{j=1}^{\theta_0}\,|\nabla^{j}\bH|\,|\nabla^{\theta_0-j+1}\bp|\,\in\bigcap_{p<\infty}L^p\:,
\ees
and thus
\be\label{fessen2}
\nabla^\perp\bp\wedge\nabla\bH\;\in\;\bigcap_{p<\infty}W^{\theta_0-1,p}\:.
\ee
Recall next the equation (\ref{prof444}) satisfied by the Gauss map, namely
\be\label{fessen1}
\Delta\bn\;=\;2\star\big(\nabla^\perp\bp\wedge\nabla\bH\big)\,+\,2\,\text{e}^{2\la}K\,\bn\,-\,2\star\text{e}^{2\la}\vec{h}_{12}\wedge(\vec{h}_{11}-\vec{h}_{22})\:.
\ee
As previously noticed, $\text{e}^{\la}\bh_{ij}$ inherits the regularity of $|\nabla\bn|$, so that (\ref{itemn1}) shows
\bes
\text{e}^{2\la}K\,\bn\,-\,\star\,\text{e}^{2\la}\vec{h}_{12}\wedge(\vec{h}_{11}-\vec{h}_{22})\;\in\bigcap_{p<\infty}W^{\theta_0-1,p}\:.
\ees
Introducing (\ref{fessen2}) and the latter into (\ref{fessen1}) now shows that 
\be\label{itemn3}
\nabla^{\theta_0+1}\bn\;\in BMO\:.
\ee

We have also seen in (\ref{item37}) that when both residues vanish, there holds
\bes
\bp\,\in\bigcap_{p<\infty}W^{\theta_0+2,p}(\di)\:.
\ees

We next study separately the cases $\theta_0<\mu+2$ and $\theta_0=\mu+2$, to see how these results might be improved.

\paragraph{\underline{Case $\theta_0<\mu+2$\,.}} Setting $\vec{\gamma}_0=\vec{0}$ in (\ref{modifres}) implies that the first residue $\vec{\beta}_0$ must vanish. In addition, since $a=0$, (\ref{item77}) shows that $\bH$ lies in $W^{1,p}(\di)$ for all $p<\infty$, thereby allowing to obtain from (\ref{residudu}) and (\ref{confdiv}) that
\be\label{confdiv3}
\Re\Big[\pz\Big(\pzb\bH\,-\,3\,\pro\pzb\bH\,+\,i\star\big(\pzb\bn\wedge\bH\big)\,-\,\text{e}^{-2\la}f\,\pz\bp\Big)\Big]\;=\;\vec{0}\quad\text{on}\:\:\di\:,
\ee
Combining (\ref{unc1})-(\ref{unc2}) yields without much effort that
\be\label{dog1}
\pi_T\pzb\bH\;=\;-\,2\,|\bH|^2\pzb\bp\,+\,i\star(\pzb\bn\wedge\bH)\:.
\ee
Hence, (\ref{confdiv3}) yields that there holds on the whole unit disk:
\bes
\dfrac{1}{2}\,\Delta\bH\;=\;\Re\Big[\pz\Big(-6\,|\bH|^2\pzb\bp\,+\,4\,i\star(\pzb\bn\wedge\bH)- \text{e}^{-2\la}f\,\pz\bp   \Big)\Big]\:.
\ees
Equivalently, using (\ref{defJ}), we may write
\be\label{booty}
\dfrac{1}{4}\,\Delta\big(2\bH+\Re(\bF_\mu)\big)\;=\;\Re\Big[\pz\Big(-6\,|\bH|^2\pzb\bp\,+\,4\,i\star(\pzb\bn\wedge\bH)-\vec{J}\Big)\Big]\:.
\ee
Since $\theta_0\le\mu+1$ and
\be\label{booty2}
\bJ\;:=\;\overline{z}^{\,\mu+1-\theta_0}\Big[z^{1-\theta_0}\text{e}^{-2u}\pz\bp-\big(z^{1-\theta_0}\text{e}^{-2u}\pz\bp\big)(0)   \Big]\:,
\ee
we see that $\bJ$ inherits the regularity of $\,z^{1-\theta_0}\text{e}^{-2u}\pz\bp$.\\

Suppose that for some integer $k\ge0$, both $\nabla\bn$ and $\bH$ lie in $C^k(\di)$. Let   $u$ be the function defined in (\ref{mulsve10}). As we have previously seen, $\Delta u$ inherits the regularity of $|\nabla\bn|^2$. Whence, 
\be\label{bts1}
u\,\in\,C^{k+2}(\di)\:.
\ee 
In particular, $\text{e}^{2u}\bH\in C^k$. Since $\text{e}^{\la}=|z|^{\theta_0-1}\text{e}^{u}$, we may write
\bes
\Delta\bp\;\equiv\;2\,\text{e}^{2\la}\bH\;=\;(\bT_k+\bR_k)\,|z|^{2(\theta_0-1)}\:,
\ees
where $\bT_k$ is the $k^\text{th}$-order Taylor polynomial of $2\text{e}^{2u}\bH$ at the origin, while $\bR_k$ is the corresponding remainder. Hence $\bp=\bp_0+\bp_1$, where
\bes
\left\{\begin{array}{lcl}\Delta\bp_0&=&|z|^{2(\theta_0-1)}\bT_k\\[1.5ex]
\Delta\bp_1&=&|z|^{2(\theta_0-1)}\bR_k\:.\end{array}\right.
\ees
Writing $\,\bT_k(x_1,x_2)\equiv\bT_k(z,\bar{z})$, there holds
\be\label{yaki0}
z^{1-\theta_0}\pz\bp_0(z,\bar{z})\;=\;\dfrac{1}{4}\int\overline{z}^{\,\theta_0-1}\bT_k(z,\bar{z})\,d\bar{z}\,+\,\vec{p}(z)\:,
\ee
where $\vec{p}$ is any meromorphic function. \\
To study the asymptotic behavior of $\bp_1$, we note first that 
\be\label{yaka7}
\nabla^j\bR_k\;=\;\text{o}(|z|^{k-j})\qquad\forall\:\:j\in\{0,\ldots,k\}\:.
\ee
We next differentiate the equation for $\bp_1$ in all variables $k$ times. With (\ref{yaka7}), one easily sees that Proposition \ref{CZpondere} is applicable and yields the representation
\be\label{yaki1}
\pz\bp_1(z,\bar{z})\;=\;\vec{P}(z)+\,\vec{\xi}(z,\bar{z})\:,
\ee
where $\vec{P}$ is a holomorphic $\C^m$-valued polynomial of degree at most $(2\theta_0+k-1)$, while the remainder $\vec{\xi}$ satisfies
\be\label{yaki2}
\nabla^{j}\vec{\xi}\;=\;\text{O}(|z|^{2\theta_0+k-j-\epsilon})\qquad\forall\:j\in\{0,\ldots,k+1\}\:\:,\:\:\forall\:\epsilon>0\:.
\ee
Combining (\ref{yaki0}) and (\ref{yaki1}) gives
\bes
z^{1-\theta_0}\pz\bp(z,\bar{z})\;=\;\vec{p}(z)+z^{1-\theta_0}\vec{P}(z)\,+\,\dfrac{1}{4}\int\overline{z}^{\,\theta_0-1}\bT_k(z,\bar{z})\,d\bar{z}\,+\,z^{1-\theta_0}\vec{\xi}(z,\bar{z})\:.
\ees
The left-hand side is bounded by assumption. The third summand on the right-hand side clearly smooth. We have also seen that the last summand on the right-hand side is bounded. Accordingly, the meromorphic function $\vec{p}+z^{1-\theta_0}\vec{P}$ is bounded and thus smooth. Altogether, we infer that the regularity of $z^{1-\theta_0}\pz\bp$ is the same as that of $z^{1-\theta_0}\vec{\xi}$. Namely, owing to (\ref{yaki2}), $\,z^{1-\theta_0}\pz\bp\in C^k(\di)$. Per (\ref{bts1}), (\ref{booty2}), and our previous discussion, there whence holds that
\be\label{regJ}
\vec{J}\,\in\,C^k(\di)\:.
\ee
Furthermore, the above arguments show that
\be\label{regphi}
\bp\,\in\,C^{k+2}(\di)\:.
\ee
As $\nabla\bn$ and $\bH$ are $k$-times continuously differentiable by hypothesis, putting (\ref{regJ}) and (\ref{regphi}) into (\ref{booty}) yields
\bes
\bH+\Re(\bF_\mu)\,\in\,C^{k+1}(\di)\:.
\ees
Returning to (\ref{defmu}), when $\,\theta_0<\mu+2$, the function $\bF_\mu$ is smooth. It then follows that $\bH$ lies in $C^{k+1}$. Introducing this fact along with (\ref{regphi}), and the hypothesis that $\nabla\bn\in C^k$ into the equation (\ref{fessen1}) for the Gauss map yields that $\nabla\bn\,\in\,C^{k+1}$. The regularities of $\bH$ and $\nabla\bn$ have both increased, and a bootstrap procedure ensues until all three functions $\bp$, $\bn$, and $\bH$ are smooth. The first step of the procedure is ensured by the fact that when $\vec{\gamma}_0=\vec{0}$ and $a=0$, the mean curvature vector lies in $\bigcap_{p<\infty}W^{1,p}\subset C^0$, while (\ref{itemn3}) guarantees that $\nabla\bn\in \bigcap_{p<\infty}W^{\theta_0,p}\subset C^0$ for any $\theta_0\ge1$.

\paragraph{\underline{Case $\theta_0=\mu+2$\,.}} 
Setting $\vec{\gamma}_0=\vec{0}$ and $a=0$ in (\ref{item77}) shows that $\bH$ lies in $W^{1,p}(\di)$ for all $p<\infty$. Moreover, (\ref{itemn3}) guarantees that $\nabla^2\bn$ lies in BMO. As in the previous paragraph, we let $u$ be the function defined in (\ref{mulsve10}). Using the local expansion (\ref{item3}), we obtain
\be\label{jbe}
\text{e}^{-2\la}f\,\pz\bp\;=\;\dfrac{\theta_0}{2}\,\text{e}^{-2u(0)}a_{\theta_0-2}\bA\;\overline{z}^{\,-1}+\,\text{O}(1)\,\in\,L^{2,\infty}\:.
\ee
On the other hand, since $f$ is anti-holomorphic, there holds as in (\ref{pety}):
\be\label{petty}
\pz\big(\text{e}^{-2\la}f\,\pz\bp\big)\;=\;\dfrac{1}{2}\,f\bH_0\qquad\text{on}\:\:\di\setminus\{0\}\:.
\ee
Because $\text{e}^{\la}\bH_0$, which is controlled by $|\nabla\bn|$, is bounded on the unit disk, we see that
\bes
|f\bH_0|\;\lesssim\;|z|^{\mu+1-\theta_0}\;=\;|z|^{-1}\,\in\,L^{2,\infty}\:.
\ees
The latter along with (\ref{jbe}) brought into (\ref{petty}) yield the identity
\be\label{petty2}
\pz\big(\text{e}^{-2\la}f\,\pz\bp\big)\;=\;\pi\dfrac{\theta_0}{2}\,\text{e}^{-2u(0)}a_{\theta_0-2}\bA\,\delta_0+\dfrac{1}{2}\,f\bH_0\qquad\text{on}\:\:\di\:.
\ee
For analogous reasons, because $\bH$ and $\nabla\bn$ belong to $W^{1,p}$ for all $p<\infty$, while $\,\text{e}^{-2\la}f\,\pz\bp$ lies in $L^{2,\infty}$, the identity (\ref{confdiv}) and the definition of (\ref{residudu}) shows that there holds on the whole unit disk
\bes
\Re\Big[\pz\Big(\pzb\bH\,-\,3\,\pro\pzb\bH\,+\,i\star\big(\pzb\bn\wedge\bH\big)\,-\,\text{e}^{-2\la}f\,\pz\bp\Big)\Big]\;=\;\pi\vec{\beta}_0\,\delta_0\:.
\ees
Combined to (\ref{petty2}), this last equation gives
\bes
\Re\Big[\pz\Big(\pzb\bH\,-\,3\,\pro\pzb\bH\,+\,i\star\big(\pzb\bn\wedge\bH\big)\Big)\Big]\;=\;\dfrac{1}{2}\,\Re(f\bH_0)\qquad\text{on}\:\:\di\:,
\ees
where we have used that
\bes
\vec{\beta}_0+\,\dfrac{\theta_0}{2}\,\text{e}^{-2u(0)}\Re\big(a_{\theta_0-2}\bA\big)\;\equiv\;\vec{\gamma}_0\;=\;\vec{0}\:.
\ees
With the help of (\ref{dog1}), we can recast the latter in the form
\bes
\dfrac{1}{2}\,\Delta\bH\;=\;\Re\Big[\pz\Big(-6\,|\bH|^2\pzb\bp\,+4\,i\star(\pzb\bn\wedge\bH)\Big)-f\bH_0\Big]\:.
\ees
As $\bH$, $\nabla\bn$, and $\nabla\bp$ all lie in $\bigcap_{p<\infty}W^{1,p}$, while $f\bH_0$ belongs to $L^{2,\infty}$, we hence infer that
\bes
\bH\,\in\,W^{2,(2,\infty)}(\di)\:.
\ees
In the special case when $\theta_0=1$, using the same techniques as those previously encountered, one derives the improvement
\bes
\bp\,\in\,W^{4,(2,\infty)}(\di)\qquad\text{and}\qquad \bn\,\in\,W^{3,(2,\infty)}(\di)\:.
\ees

\subsection{Special cases}

\subsubsection{Regular points (no branch)}

When the origin is a regular point, the modulus of the conformal parameter $|\la|$ is bounded through the origin (exactly as when the origin is a branch point of order one, i.e. $\theta_0=1$, but unlike in the case when the origin is a branch point of order $\theta_0\ge2$, for which $\lim_{x\rightarrow0}\la(x)=-\infty$).\\
We have seen in (\ref{item37}) that $\bH$ lies in $BMO(\di)$ when $\theta_0=1$. Moreover, when $\theta_0=1$, then $a$ is automatically zero, so that the second residue $\vec{\gamma}$ at the origin vanishes. We now verify that the modified first residue $\vec{\gamma}_0$ must also vanish. When $\theta_0=1$, we have seen that
\bes
\pz\bp\;=\;\dfrac{1}{2}\,\bA+\text{O}(|z|)\qquad\text{and}\qquad \text{e}^{-2\la}\;=\;\text{e}^{-2u(0)}+O(|z|)\:.
\ees
Since the anti-holomorphic multiplier function satisfies
\bes
f(\bar{z})\;=\;\dfrac{a_{-1}}{\bar{z}}+f_0(\bar{z})\qquad\text{with}\quad a_{-1}\in\C\quad\text{and}\quad f_0\in C^\infty(\di)\:,
\ees
it follows that
\begin{eqnarray*}
\dfrac{1}{2}\,f\bH_0\;=\;\pz\big(\text{e}^{-2\la}f\pz\bp\big)\,-\,a_{-1\,}\text{e}^{-2\la}\pz\bp\,\pz\dfrac{1}{\bar{z}}\qquad\text{on}\:\:\di\:,
\end{eqnarray*}
so that
\be\label{befroi}
\int_{\di}\Re\big(f\bH_0\big)\;=\;-\,\pi\,\text{e}^{-2u(0)}\Re\big(a_{-1}\bA\big)\,+\,\dfrac{1}{2}\int_{\partial\di}\vec{\nu}\cdot\big(\text{e}^{-2\la}M_f\nabla^\perp\bp\big)\:,
\ee
with
\bes
M_f\;:=\;\left(\begin{array}{rc}-\Im(f)&\Re(f)\\[.6ex]\Re(f)&\Im(f)\end{array}\right)\:.
\ees
If the origin is a regular point, then by definition the constrained Willmore equation holds on the whole unit disk:
\bes
\text{div}\,\Big[\nabla\bH\,-\,3\,\pro\nabla\bH\,+\,\star\,\big(\nabla^\perp\bn\wedge\bH\big)\Big]\;=\;2\,\Re\big(\bH_0f\big)\:.
\ees
Upon integrating over $\di$, using the definition of the first residue (cf. (\ref{residudu})), and introducing (\ref{befroi}), we find
\bes
\vec{\beta}_0\;=\;-\,\dfrac{1}{2}\,\text{e}^{-2u(0)}\Re\big(a_{-1}\bA\big)\:,
\ees
so that indeed the modified first residue $\vec{\gamma}_0$ defined in (\ref{modifres}) vanishes when $\theta_0=1$ and the origin is a regular point. \\

Knowing now that both residues $\vec{\gamma}_0$ and $\vec{\gamma}$ vanish, we call upon the results of Section \ref{wbrv} to infer that
\begin{itemize}
\item if $\mu\ge0$ (i.e. $f$ is regular at the origin), then $\bp$ is smooth ;
\item if $\mu=-1$ (i.e. $f$ is singular at the origin), then there holds
\bes
\bp\,\in\,W^{4,(2,\infty)}\quad,\quad \bn\,\in\,W^{3,(2,\infty)}\quad,\quad \bH\,\in\,W^{2,(2,\infty)}\:,
\ees
and thus in particular that $\bp$ is $C^{2,\alpha}(\di)$ for all $\alpha\in[0,1)$.
\end{itemize}

\subsubsection{Surfaces of specific types}
As we have seen in the introduction, constrained Willmore surfaces englobe various types of commonly studied surfaces. They include of course Willmore surfaces (comprising {\it inter alia} minimal surfaces). The Willmore equation is obtained from the constrained Willmore equation by setting the multiplier function $f\equiv0$. In that case, one unsurprisingly checks that our Theorems \ref{Th4} and \ref{Th44} are the main results of \cite{BR3}.\\

The other important subclass of constrained Willmore surfaces are the {\it parallel mean curvature surfaces} (generalizing constant mean curvature surfaces to higher codimension). These surfaces have the property that
\be\label{pmccond}
\pro\nabla\bH\;\equiv\;\vec{0}\:.
\ee
As shown in Lemma \ref{PMC}, they satisfy the conformal Willmore equation with multiplier $\,f=-2\,\text{e}^{2\la}\bH\cdot\bH_0^*$. \\
The residues associated with such a surfaces necessarily vanish. Indeed, using the identity (\ref{unc1}) and Proposition \ref{Th1}, we infer the following estimate for the tangential projection:
\bes
\big|\pi_T\nabla\bH\big|\;\lesssim\;|\nabla\bn||\bH|\;\lesssim\;\text{e}^{-\la}|\nabla\bn|^2\;=\;\text{O}(|z|^{1-\theta_0-\epsilon})\qquad\forall\:\epsilon>0\:.
\ees
Combining the latter to (\ref{pmccond}) implies that 
\bes
\left\{\begin{array}{lcl}
\bH\in\,\bigcap_{p<\infty}W^{1,p}&&\text{if}\:\:\:\theta_0=1\\[1ex]
\bH\,=\,\text{O}(|z|^{2-\theta_0-\epsilon})&&\text{if}\:\:\:\theta_0\ge2\:.
\end{array}\right.
\ees
In the second case, we can repeat exactly the same routine, with $\theta_0-1$ in place of $\theta_0$. The procedure continues finitely many times, until reaching that in all cases, there holds $\,\bH\in\bigcap_{p<\infty}W^{1,p}$. As seen in (\ref{item77}), this guarantees that both residues associated to the surface at the origin necessarily vanish. On the other hand, we obtain from the boundedness of $\bH$ and the fact that $\text{e}^{\la}|\bH_0|\lesssim|\nabla\bn|=\text{O}(|z|^{-\epsilon})$ for all $\epsilon>0$, that
\bes
|f(\bar{z})|\;\lesssim\;|z|^{\theta_0-1-\epsilon}\qquad\forall\:\epsilon>0\:.
\ees 
As $f$ is anti-holomorphic, its order at the origin thus satisfies $\mu\ge\theta_0-1$. We may now call upon Theorem \ref{Th44} to conclude that a parallel mean curvature surface is smooth throughout its interior branch point.

\bigskip

\renewcommand{\theequation}{A.\arabic{equation}}
\renewcommand{\theTh}{A.\arabic{Th}}
\renewcommand{\theProp}{A.\arabic{Prop}}
\renewcommand{\theLm}{A.\arabic{Lm}}
\renewcommand{\theCo}{A.\arabic{Co}}
\renewcommand{\theRm}{A.\arabic{Rm}}
\renewcommand{\theequation}{A.\arabic{equation}}
\setcounter{equation}{0} 
\reset
\appendix
\section{Appendix}
\subsection{Notational Conventions}

We append an arrow to all the elements belonging to $\R^m$. To simplify the notation, by $\bp\in X(\di)$ is meant $\bp\in X(\di\!,\R^m)$ whenever $X$ is a function space. Similarly, we write $\nabla\bp\in X(\di)$ for $\nabla\bP\in \mathbb{R}^2\otimes X(\di\!,\R^m)$.\\[1.5ex]
Although this custom may seem at first odd, we allow the differential operators classically acting on scalars to act on elements of $\R^m$. Thus, for example, $\nabla\bp$ is the element of $\R^2\otimes\R^m$ that can be written $(\px\bp,\py\bp)$. If $S$ is a scalar and $\bR$ an element of $\R^m$, then we let
\begin{eqnarray*}
\bR\cdot\nabla\bP&:=&\big(\bR\cdot\px\bP\,,\,\bR\cdot\py\bP\big)\:\\[1ex]
\nabla^\perp S\cdot\nabla\bP&:=&\px S\,\py\bp\,-\,\py S\,\px\bp\:\\[1ex]
\nabla^\perp\bR\cdot\nabla\bP&:=&\px\bR\cdot\py\bp\,-\,\py\bR\cdot\px\bp\:\\[1ex]
\nabla^\perp\bR\wedge\nabla\bP&:=&\px\bR\wedge\py\bp\,-\,\py\bR\wedge\px\bp\:.
\end{eqnarray*}
Analogous quantities are defined according to the same logic. \\

Two operations between multivectors are useful. The {\it interior multiplication} $\res$ maps a pair comprising a $q$-vector $\gamma$ and a $p$-vector $\beta$ to a $(q-p)$-vector. It is defined via
\bes
\langle \gamma\res\beta\,,\alpha\rangle\;=\;\langle \gamma\,,\beta\wedge\alpha\rangle\:\qquad\text{for each $(q-p)$-vector $\alpha$.}
\ees
Let $\al$ be a $k$-vector. The {\it first-order contraction} operation $\bul$ is defined inductively through 
\bes
\alpha\bul\beta\;=\;\alpha\res\beta\:\:\qquad\text{when $\beta$ is a 1-vector}\:,
\ees
and
\bes
\alpha\bul(\beta\wedge\gamma)\;=\;(\alpha\bul\beta)\wedge\gamma\,+\,(-1)^{pq}\,(\alpha\bul\gamma)\wedge\beta\:,
\ees
when $\beta$ and $\gamma$ are respectively a $p$-vector and a $q$-vector.

\subsection{Miscellaneous Facts}

\begin{Lm}\label{PMC}
Parallel mean curvature immersions are constrained Willmore with multiplier $f=-2\text{e}^{2\la}\bH\cdot\bH_0^*$.
\end{Lm}
$\textbf{Proof.}$ The constrained Willmore equation (\ref{wil1}) reads
\be\label{wil11}
\Delta_\perp\bH\,+\,2\,\Re\big((\bH\cdot\bH^*_0)\bH_0\big)\:=\:\text{e}^{-2\la}\Re({\bH_0}f)\:,
\ee
with
\bes
-\,\Delta_\perp\bH\::=\:\text{e}^{-2\la}\,\pi_{\bn}\,\text{div}\,\pi_{\bn}\nabla\bH\:.
\ees
By definition, a parallel mean curvature immersion has $\pro\nabla\bH\equiv\vec{0}$. Whence (\ref{wil11}) is satisfied with $f:=2\text{e}^{2\la}\bH\cdot\bH_0^*$. Remains thus to verify that $f$ is indeed anti-holomorphic. To this end, we recall the Codazzi equation written in the form\footnote{cf. Lemma A.3 in \cite{BR1}. {\it Caveat:} what was known as $\bH_0$ in \cite{BR1} is known as $\bH_0^*$ in the present work.}:
\bes
\text{e}^{-2\la}\pzb\big(\text{e}^{2\la}\bH\cdot\bH_0^*\big)\;=\;\bH\cdot\pz\bH\,+\,\bH_0\cdot\pzb\bH\:.
\ees
As $\bH$ and $\bH_0$ are normal vectors, while $\pz\bH$ and $\pzb\bH$ are tangent vectors by assumption, the latter confirms, as desired, that $f$ is anti-holomorphic.\\[-1.5ex]

\hfill $\blacksquare$

\begin{Lm}\label{GW}
Let $u\in C^2(\di\setminus\{0\})$ solve the equation
\be\label{vdb}
\text{div}\,\big(\nabla u(x)+V(x,u)+T(x)\big)\;=\;0\qquad\text{on}\:\: \di\setminus\{0\}\:.
\ee
Assume that for some integer $a\ge1$, there holds
\bes
|x|^aT\:\:,\:\:|x|^aV\:\:,\:\:|x|^{a-1}u\:\:\in\, L^2(\di)\:.
\ees
Then we have
\bes
|x|^{a}\nabla u\,\in\,L^2(\di)\:.
\ees
\end{Lm}
$\textbf{Proof.}$
We fix a point $y_0\neq0$, set $d:=\frac{1}{2}|y_0|$, and we define a $C^2$ cut-off function $\eta:=|x|^{a}\mu$, where $\mu$ is the standard smooth cut-off function
\be\label{est1}
\mu\,=\,\Bigg\{\begin{array}{ll}1\:,&D_{d/2}(y_0)\\[1.5ex]0\:,&D^c_d(y_0)\end{array},\quad\text{and}\qquad \Bigg\{\begin{array}{l}|\nabla\mu|\,\lesssim\,d^{-1}\\[1.5ex]|\nabla^2\mu|\,\lesssim \,d^{-2}\end{array}\quad\text{on}\:\:\:D_d(y_0)\:.
\ee

\noindent
We will need the function 
\bes
F_y(x)\;:=\;\dfrac{1}{4\pi}\log|x-y|\:.
\ees
Applying the test-function $\,\eta F_y$\, to the equation (\ref{vdb}) yields after performing a few elementary manipulations and integrating by parts, the identity
\bes
|y|^{a}u(y)\;=\;-\,\int_{D_d(y)}u\,\big(2\,\nabla_x F_y\cdot\nabla_x\eta+F_y\Delta_x\eta\big)\;+\;\int_{D_d(y)}(V+T)\cdot\nabla_x(\eta F_y)\:.
\ees
We next apply the operator $\nabla_y$ to both sides of the latter to obtain
\begin{eqnarray}\label{murdoch}
\nabla\big(|y|^{a}u(y)\big)&=&\int_{D_d}\nabla_{xy}^2F_y\big((V+T)\eta-2u\nabla_x\eta\big)+\nabla_y F_y\big((V+T)\cdot\nabla_x\eta-u\Delta_x\eta \big)\nonumber\\[1ex]
&=&\Omega*\Big[\big((V+T)\eta-2u\nabla_x\eta\big)\chi_{D_d(y)}(x)\Big]\nonumber\\[-0ex]
&&\hspace{.75cm}+\:\:\dfrac{1}{4\pi}\int_{D_d\setminus D_{d/2}}\big((V+T)\cdot\nabla_x\eta-u\Delta_x\eta \big)\dfrac{x-y}{|x-y|^2}\nonumber\\[1ex]
&=:&I_1(y)+I_2(y)\:,
\end{eqnarray}
where $\Omega$ is the standard Calderon-Zygmund kernel:
\bes
\Omega(z)\;:=\;\dfrac{|z|^2\,\mathbb{I}_2\,-\,2\,z\otimes z}{4\pi|z|^4}\:.
\ees
On the annulus $D_d\setminus D_{d/2}$, there holds
\bes
\dfrac{d}{2}<|x-y_0|<d\qquad\text{and}\qquad d<|x|<3\,d\:.
\ees
Hence, on the annulus,
\bes
|\nabla\eta|(x)\;\lesssim\;|x|^{a-1}+d^{-1}|x|^{a}\;\lesssim\;d^{-1}|x|^{a}
\qquad\text{and}\qquad
|\Delta\eta|(x)\;\lesssim\;d^{-1}|x|^{a-1}\:.
\ees
Accordingly, the second integral in (\ref{murdoch}) is estimated (up to an irrelevant constant) by
\bes
|I_2|(y_0)\:\lesssim\:M\big[|x|^{a}V\big](y_0)\,+\,M\big[|x|^{a}T\big](y_0)\,+\,M\big[|x|^{a-1}u\big](y_0)\:,
\ees
where $M$ denotes the standard maximal function, and we have used a classical estimate which bounds the Riesz transform by the maximal function (\cite{Zi}, Lemma 2.8.3). Accordingly, we find
\begin{eqnarray*}
&&\hspace{-1cm}\Vert I_2\Vert_{L^p(D_{2/3}(0))}\nonumber\\[1ex]
&\lesssim&\big\Vert M\big[|x|^{a}V\big]\big\Vert_{L^p(D_{2/3}(0))}+\big\Vert M\big[|x|^{a}T\big]\big\Vert_{L^p(D_{2/3}(0))}+\big\Vert M\big[|x|^{a-1}u\big]\big\Vert_{L^p(D_{2/3}(0))}\nonumber\\[1ex]
&\lesssim&\big\Vert |x|^{a}V\big\Vert_{L^p(D_{2/3}(0))}+\big\Vert|x|^{a}T\big\Vert_{L^p(D_{2/3}(0))}+\big\Vert |x|^{a-1}u\big\Vert_{L^p(D_{2/3}(0))}\:,
\end{eqnarray*}
for any $p\in(1,\infty)$. On the other hand, standard $L^p$ estimates for the convolution with the Calderon-Zygmund kernel yield
\bes
\Vert I_1\Vert_{L^p(D_{2/3}(0))}\:\lesssim\:\big\Vert|x|^{a-1}u\big\Vert_ {L^p(D_{2/3}(0))}\:.
\ees
Combining the last two estimates thus gives
\bes
\Big\Vert \nabla\big(|x|^{a}u\big)\Big\Vert_{L^p(D_{2/3}(0))}\;\lesssim\;\big\Vert|x|^{a}V\big\Vert_{L^p(\di)}+\,\big\Vert|x|^{a}T\big\Vert_{L^p(\di)}+\,\big\Vert|x|^{a-1}u\big\Vert_{L^p(\di)}\:,
\ees
from which it easily follows that
\bes
\big\Vert|x|^{a}\nabla u\big\Vert_{L^p(D_{2/3}(0)))}\;\lesssim\;\big\Vert|x|^{a}V\big\Vert_{L^p(\di)}+\,\big\Vert|x|^{a}T\big\Vert_{L^p(\di)}+\,\big\Vert|x|^{a-1}u\big\Vert_{L^p(\di)}\:.
\ees
In particular, using $p=2$ and the given hypotheses, we infer the announced statement.\\[-3ex]

\hfill $\blacksquare$\\

The next two results are proved in detail in the Appendix of \cite{BR3}. We content ourselves here with stating them only.

\begin{Prop}\label{CZpondere}
Let $u\in C^2(\di\setminus\{0\})$ solve
\bes
\Delta u(x)\;=\;\mu(x)f(x)\qquad\text{in\:\:\:} \di\:,
\ees
where $f\in L^p(\di)$ for some $p>2$. The weight $\mu$ satisfies
\bes
|\mu(x)|\;\simeq\;|x|^{b}\qquad\quad\text{for some}\:\:\: b\in\mathbb{N}\:.
\ees
Then 
\begin{itemize}
\item[(i)] there holds\footnote{$\overline{x}$ is the complex conjugate of $x$. We parametrize $\di$ by $x=x_1+i\,x_2$, and then $\overline{x}:=x_1-i\,x_2$. In this notation, $\nabla u$ in (\ref{stmt}) is understood as $\partial_{x_1}u+i\,\partial_{x_2}u$.}
\be\label{stmt}
\nabla u(x)\;=\;P(\overline{x})\,+\,|\mu(x)|\,T(x)\:,
\ee
where $P(\overline{x})$ is a complex-valued polynomial of degree at most $b$, and near the origin $\,T(x)=\text{O}\big(|x|^{1-\frac{2}{p}-\epsilon}\big)$ for every $\epsilon>0$. \\[-1.5ex]
\item[(ii)] furthermore, if $\,\mu\in C^1(\di\setminus\{0\})$ and if
\be\label{hypw2}
|x|^{1-b}\,\nabla\mu(x)\,\in\,L^{\infty}(\di)\:,
\ee
there holds
\bes\label{stmt2}
\nabla^2 u(x)\;=\;\nabla P(\overline{x})\,+\,|\mu(x)|\,Q(x)\:,
\ees
where $P$ is as in (i), and 
\bes
Q\;\in\;L^{p-\epsilon}(\di,\C^2)\qquad\quad\forall\:\:\epsilon>0\:.
\ees
As a $(2\times2)$ real-valued matrix, $Q$ satisfies in addition
\bes
\text{Tr}\;Q\;\in\;L^p(\di)\:.
\ees
\end{itemize}
\end{Prop}

\smallskip

\begin{Prop}\label{CZcoro}
Let $u\in C^2(\di\setminus\{0\})$ solve
\bes
\Delta u(x)\;=\;\mu(x) f(x)\qquad\text{in\:\:\:} \di\:,
\ees
where
\bes
|f(x)|\;\lesssim\;|x|^{n+r}\qquad\text{and}\qquad|\mu(x)|\;\simeq\;|x|^{b}\:,
\ees
for two non-negative integers $n$ and $b$ ; and for some $r\in(0,1)$. \\
Then
\bes
\nabla u(x)\;=\;P(\overline{x})\,+\,|\mu(x)|\,T(x)\:,
\ees
where $P$ is a complex-valued polynomial of degree at most $(b+n+1)$, and near the origin $T(x)=\text{O}(|x|^{n+1+r-\epsilon})$ for every $\epsilon>0$.\\[1ex]
If in addition $\mu$ satisfies (\ref{hypw2}), then $\,|x|^{-(n+r)}|\mu|^{-1}\nabla\big(|\mu|T\big)$ belongs to $L^p$ for all finite $p$. Furthermore, there holds the estimate
\bes
\big|\text{Tr}\;\nabla\big(|\mu(x)|T(x)\big)\big|\;\;\lesssim\;\;|x|^{n+r}|\mu(x)|\:.
\ees
\end{Prop}

\eject

\end{document}